\documentclass[reqno,11pt]{amsart}         

\usepackage{graphicx}        
\usepackage{amsthm,amssymb,amsfonts,amsmath}

\usepackage{bbm,color}
\usepackage{hyperref}     
         
\newtheorem{theo}{Theorem}[section]

\newtheorem{rem}{Remark}[section]
\newtheorem{lem}{Lemma}[section]
\newtheorem{cor}{Corollary}[section]

\newtheorem{ex}{Example}[section]

\newcommand{\rme}{{\rm e}}
\newcommand{\rmd}{{\rm d}}
\newcommand{\la}{{\lambda}}
\newcommand{\ep}{{\varepsilon}}
\newcommand{\1}{\mathbf{1}}

\newcommand{\sF}{{\mathcal F}}

\newcommand{\ttT}{{\mathtt T}}
\newcommand{\PP}{{\mathbb P}}

\newcommand{\RR}{{\mathbb R}}

\newcommand{\EE}{{\mathbb E}}

\newcommand{\bp}{{\mathbf p}}
\newcommand{\emm}{{ \mathit m}}
\newcommand{\Ell}{{ \mathtt L}}

\begin{document}

\title{Kac-L\'evy processes}\thanks{To appear in \emph{Journal of Theoretical Probability}}
\author{Nikita Ratanov}
 \address{N. Ratanov, Facultad de Econom\'{\i}a, Universidad del Rosario, 
Calle 12c, No.4-69, Bogot\'a,  D. C. Colombia}
\email{nikita.ratanov@urosario.edu.co}
\date{}

 \maketitle

\begin{abstract}
Markov-modulated L\'evy processes with two different regimes of restarting are studied.
These regimes correspond to the completely renewed process 
and  to the process of Markov modulation, accompanied  by jumps.
We give explicit expressions for the L\'evy-Khintchine exponent in the case of a two-state underlying Markov chain.

For the renewal case, the limit distributions (as $t\to\infty$) are obtained.
In the case of processes with jumps, 
we present some results for the exponential functional. 
\end{abstract}\vskip5mm



\section{Introduction: Kac-L\'evy process}

New classes of stochastic processes appearing 
in various fields of applications could be constructed as follows. 
First, we have the set of underlying processes, the main blocks of the construction, 
second, the switching mechanism between these blocks  should be described, 
and third, the new starting points of these underlying processes could be determined in various manners.

For instance, let $\eta_n=\eta_n^x(t)\in\mathbb X,\;t\geq0, \; n\geq0,\;$ be the sequence 
of random processes, $\{T_n\}_{n\geq0}$ be the sequence of switching times and
  $x=x_n\in  \mathbb X,\;n\geq0,$ be a  sequence of starting points.  
Here $ \mathbb X$ is the phase space (some topological vector space, in general).
The composite process $X$ is defined by
\begin{equation}\label{def:X}
X(t)=\eta_n^{x}(t-T_n)\qquad \text{for}\qquad t\in[T_n, \; T_{n+1}),\qquad n\geq0.
\end{equation}
This construction can be accompanied by a simple extension  that 
allows for additional jumps at the times of the state change.

In this paper, the sequence $\{T_n\}_{n\geq0}$  
 is defined by a finite-state Markov process 
$\ep=\ep(t)\in E,\;t\geq0,$ that is 
$T_n$ is the switching time of $\ep$ from one state to another.
Assume that 
 $\{\eta_n\}_{n\geq0}$ is a sequence of real-valued one-dimensional
  independent L\'evy processes  independent of $\ep.$

Let us summarise the most important known results.
Assume that  $\eta(t),$ $t\geq0,$ is a real-valued one-dimensional L\'evy process defined 
on the probability space $(\Omega, \sF, \PP),$ 
 that is the process satisfying the following properties:
\begin{itemize}
  \item The paths are a. s. right-continuous with left limits, and $\PP(\eta(0)=0)=1;$
  \item $\eta(t)-\eta(s)\stackrel{D}{=}\eta(t-s),\;0\leq s\leq t;$
 \item  $\eta(t)-\eta(s)$ is independent of the history $\{\eta_u~|~u\leq s\}.$
\end{itemize}
The jump part of the L\'evy process $\eta$ can be characterised by means of 
the L\'evy measure $\mathrm\Pi$ (intensity measure),
\begin{equation*}
\mathrm\Pi(A):=\EE\left[J([0, 1]\times A)\right],\qquad A\in\mathcal B(\RR\setminus\{0\}),
\end{equation*}
satisfying the condition
\begin{equation*}
\int_\RR(1\land x^2)\mathrm\Pi(\rmd x)<\infty.
\end{equation*}
Here the jump measure  $J(I, A),\;I\in\mathcal B\left([0, \infty)\right),\;A\in\mathcal B(\RR),$
is the random measure which is defined by 
\[J(I, A)=\sum_{s\in I}\1_{\Delta\eta(s)\in A\setminus \{0\}}.\]

In general, process $\eta$ is completely identified by the triplet
$\langle\gamma,\;\sigma^2,\;\mathrm\Pi\rangle,$ see e.g. \cite{kyprianou,sato}.
This means that $\eta(t)$ may be decomposed into a sum of independent terms
(L\'evy-It\^o decomposition):
\begin{equation*}
\eta(t)=\gamma t+\sigma W(t)+Y(t)+M(t),\qquad t\geq0.
\end{equation*}
Here $\gamma$ and $\sigma$ are constants,
$W(t),\;t\geq0,$ is the standard Brownian motion, $Y=Y(t)$ is the compound Poisson process,
\[
Y(t)=\int_0^t\int_{|y|\geq R}yJ(\rmd s, \rmd y)=\sum_{n=1}^{N_t}Y_n,
\]
where $\{Y_n\}$ are independent random variables with common distribution 
supported on $\{y: |y|\geq R\},$
and $M=M(t)$ is the square-integrable martingale,
\[
M(t) 
=\int_0^t\int_{|y|<R}yJ(\rmd s, \rmd y)-t\int_{|y|<R}y\mathrm\Pi(\rmd y).
\]

The distribution of $\eta(t)$ is determined by the characteristic function 
$\EE\left[\rme^{\mathrm i\theta\eta(t)}\right]$ which is given by $\exp(-t\psi(\theta)).$
According to the L\'evy-Khintchine formula, $\psi(\theta)$ can be written as 
\begin{align}\nonumber
    \psi(\theta)&  = \left\{
     - \mathrm i \gamma\theta
    +\frac12\sigma^2\theta^2\right\}\\ 
     \nonumber
    &  +\left\{
    \int_{\{|x|\geq R}(1-\rme^{\mathrm i\theta x})\mathrm\Pi(\rmd x)
    \right\}\\
    \nonumber
    &+\left\{
    \int_{\{0<|x|<R\}}(1-\rme^{\mathrm i\theta x}+\mathrm i\theta x)\mathrm\Pi(\rmd x)
    \right\},\qquad \theta\in\RR.
\end{align}
Notice that the real part of $\psi(\theta)$ is nonnegative,
\begin{equation}
\label{eq:Repsi}
\mathrm{Re}[\psi(\theta)]
=\frac12\sigma^2\theta^2+ \int_\RR\left(1-\cos\theta x\right)\mathrm\Pi(\rmd x)\geq0,
\qquad \forall\theta.
\end{equation}
We will also use the L\'evy-Laplace exponent
\[\ell(\xi)=\psi(\mathrm i\xi)=\EE\left[\exp\left(-\xi X(t)\right)\right];
\qquad
\EE\left[\rme^{-\xi\eta(t)}\right]=\exp\left(-t\ell(\xi)\right).
\]

A Markov-modulated L\'evy process appears 
when we try to define a composite process $X$ 
that involves switching between the L\'evy components
once the underlying Markov process $\ep$ is switched.
 It is assumed that in the interval
$t\in[T_n,\; T_{n+1})$ where $\ep$ is constant, $\ep(t)\equiv i,$
the process $X(t)$ is developing as a   L\'evy  process with the triplet 
$\langle\gamma_i,\;\sigma_i^2,\;\mathrm\Pi_i\rangle,\;i\in E,$  see \eqref{def:X}.
Aspiring  to  simplicity, 
in this paper we study the case of Markov modulation with two states, i. e. $E=\{0, 1\}.$
 
In recent decades, the Markov-modulated L\'evy processes with various extensions become
 very popular among researchers from different points of view. 
Terminology also depends on the scope of application.

This type of process is called the  \emph{Markov additive processes}, if it is used
in queuing theory, see \cite{asmussen,EJP2014}. Regarding the financial market modeling, this 
process becomes a \emph{L\'evy regime-switching model},              see \cite{boyar1,chev,Chou,Elliott2,Hainaut2,Hainaut1,Elliott1},
or a \emph{Markov-switching model}, \cite[Section 9.5]{Epps}. 
On the other hand, Markov modulated piecewise linear processes 
have been extensively studied, starting with the seminal paper by Marc Kac, \cite{Kac}.

To unify the terminology and to avoid misunderstandings, 
I propose the term \emph{Kac-L\'evy process} which can be applied 
to the composite process $X,$ \eqref{def:X}, with  L\'evy underlying blocks.

In this paper, the two different extensions are studied in detail. 

First, in Section \ref{sec:renewal} we assume that the Kac-L\'evy process is completely
renewed after each switch, that is, the process is restarted from a new random point.
In this case \eqref{def:X} becomes
\begin{equation}
\label{eq2:X}
X(t)=x_n+\eta_{n}(t-T_n),\qquad T_n\leq t<T_{n+1},\qquad n\geq0,
\end{equation}
where the starting points $\{x_n\}$ are independent random variables, 
independent of the L\'evy processes
$\{\eta_{n}(t)\},$ $t\geq0.$
In this case, the composite  process $X$
 resembles well-studied Markovian growth-collapse processes, 
which presume 
 a constant trend with additive or multiplicative down jumps.  see e.g. \cite{Boxma}. 
 Such models occur
in insurance mathematics and related fields, see  \cite[XIV-5]{asmussen} 
 or \cite[Chapters 5 and 11]{schmidt}, and in
production/inventory models studied by Shanthikumar and Sumita, \cite{shanthikumar},
 among others.

Second, in Section \ref{sec:jump} we consider the Markov modulation accompanied by jumps.
We study a Kac-L\'evy process $X$ extended with a jump component,
\begin{equation*}
\begin{aligned}
X(t)=&X(T_n)+Y_{n-1}+\eta_n(t-T_n),\qquad &T_n\leq t<T_{n+1},&\qquad n\geq1,\\
X(t)=&\eta_0(t),\qquad &0\leq t<T_1,&
\end{aligned}
\end{equation*}
where the underlying L\'evy processes $\{\eta_n(t),\;t\geq0\}_{n\geq0}$ 
and jump magnitudes $\{Y_n\}_{n\geq0}$ are independent. 
These processes are also called  \emph{Markov additive processes}, see \cite[Proposition 3.1]{EJP2014}.

In both cases we present the characteristic/moment generating functions 
of the distribution of $X(t)$.

 In particular, if   the underlying processes are determined by
$\sigma_i=0,\; \mathrm\Pi_i=0$ and $\gamma_i=c_i,$
that is $\eta_n(t)=c_it,\;i\in\{0, 1\},$ 
then we have the important  well-studied example of
the  composite
process $X(t),\;t\geq0,$ which is called (integrated)
\emph{telegraph} process (\emph{Goldstein-Kac process}),  \cite{Gold,Kac}, 
see also \cite{KR,APT}, 
\begin{equation}\label{def:telegraph}
X(t)=\int_0^tc_{\ep(s)}\rmd s.
\end{equation}
Process \eqref{def:telegraph} is also called  \emph{piecewise linear Markov process}
or \emph{Markovian fluid}, \cite{asmussen},
This is the persistent random motion switching between two trends at Markovian instants.

The transition densities $p_i(t, y)$ of the telegraph process $X(t),$ 
\[
p_i(t, y)\rmd y=\PP\{X(t)\in\rmd y~|~\ep(0)=i\},\qquad i\in\{0, 1\},
\]
follow the coupled integral equations,
\begin{equation}
\label{eq:eq-telegraph}
\left\{\begin{aligned}
p_0(t, y)=&\rme^{-\la_0t}\delta(x-c_0t)+\int_0^t\la_0\rme^{-\la_0\tau}p_1(t-\tau, y-c_0\tau)\rmd\tau,\\
p_1(t, y)=&\rme^{-\la_1t}\delta(x-c_1t)+\int_0^t\la_1\rme^{-\la_1\tau}p_0(t-\tau, y-c_1\tau)\rmd\tau,
\end{aligned}\right.\end{equation}
where $\la_0$ and $\la_1$ are the switching intensities.

Further, the moment generating functions 
$\Ell_i(t, \xi):=\EE\left[\rme^{-\xi X(t)}~|~\ep(0)=i\right],$ $i\in\{0, 1\},$
of $X(t)$ can be expressed by
\begin{equation*}
\left\{
\begin{aligned}
    \Ell_0(t, \xi)&=\rme^{-t(\la+a\xi)}\left(\cosh(tD)+(\la-c\xi)\frac{\sinh(tD)}{D}\right),   \\
    \Ell_1(t, \xi)&=\rme^{-t(\la+a\xi)}\left(\cosh(tD)+(\la+c\xi)\frac{\sinh(tD)}{D}\right),  
\end{aligned}\right.
\end{equation*}
where $a=(c_0+c_1)/2,\;c=(c_0-c_1)/2,\;\la=(\la_0+\la_1)/2,\;\mu=(\la_0-\la_1)/2$ 
and $D=D(\xi)=\sqrt{(\mu+c\xi)^2+\la_0\la_1},$
see \cite{STAPRO82}.  Occupation time distribution for such a process has been analysed by
\cite{SPA121}.
For the detailed analysis and the properties of the distribution of such a process see \cite{KR}.

In Sections \ref{sec:exponential} and \ref{sec:exponential2}, the distribution 
of the exponential functional of the Kac-L\'evy process with jumps is studied
and some explicit formulae are presented  for the Goldstein-Kac process.
Some other rare  examples with the similar explicit formulae of this type can be found in
\cite{BBY,BY,B-NMR,Kuznetsov}.
\section{Kac-L\'evy processes  with renewal starting points}\label{sec:renewal}
\setcounter{equation}{0}
\setcounter{figure}{0}

Consider a Markov-modulated L\'evy process $X=X(t)$
based on the two-state Markov process $\ep=\ep(t)\in\{0, 1\},$ $t\geq0,$
and on the independent of $\ep$  sequence $\{\eta_n(t)\}_{n\geq0},\; t\geq0,$
of independent L\'evy processes with the alternating marginal distributions 
$q^t_0(\rmd y)$ and $q_1^t(\rmd y).$ Let $\psi_0$ and $\psi_1$ be the corresponding 
L\'evy-Khintchine exponents.

At every switching time $T_n,$ process $X$  starts  from a renewed random  point $x_n$:
\begin{equation*}
X(t)=x_n+\eta_{n}(t-T_n),\qquad T_n\leq t<T_{n+1},\quad n\geq0,
\end{equation*}
see \eqref{eq2:X}. Here $x_n,\; n\geq0,$ are  independent random variables
with the alternating distributions
$g_0=g_0(\rmd x)$ and $g_1=g_1(\rmd x).$
Suppose that the
variables $\{x_n\}_{n\geq0}$  are independent of   $\{\eta_n\}_{n\geq0}$ and  $\ep.$ 

Therefore, 
\begin{equation}
\label{def:dnXt-renewal}
\PP\{X(t)\in\rmd z~|~T_n\leq t<T_{n+1}\}=g_{\ep_n}*q_{\ep_n}^{t-T_n}(\rmd z).
\qquad \ep_n=\ep(T_n).
\end{equation}
Here $*$ means convolution, such that for any test-function $\varphi$
\[
\int_{\RR^1_z}\varphi(z)\left[g*q\right](\rmd z)=\int_{\RR^2_{(x,y)}}\varphi(x+y)g(\rmd x)q(\rmd y).
\]

The distribution of $X(t)$  is determined by the density vector-function 
\[\bp^X(t, \rmd y)=(p^X_0(t,  \rmd y),\;p^X_1(t,  \rmd y))^{\tt T},\] with the entries
$p^X_i(t,  \rmd y)=\PP\{X(t)\in\rmd y~|~\ep(0)=i\},\;i\in\{0, 1\}.$  
We have the following representation of $\bp^X(t, \rmd y).$
\begin{theo}\label{theo:renewal}
Let
$\Lambda^{\mathrm{diag}}$ be the diagonal matrix of the switching intensities\textup{,}
\begin{equation*}
\Lambda^{\mathrm{diag}}=
\begin{pmatrix}
   \la_0   & \quad 0  \\ \\
    0  &  \quad \la_1
\end{pmatrix}\end{equation*}
and
\begin{equation}
\label{def:B}
\mathcal B(t)
=\frac{1}{2\la}
\begin{pmatrix}
     1-\rme^{-2\la t} &  \quad 1+\dfrac{\la_0}{\la_1}\rme^{-2\la t} \\ \\
    1+\dfrac{\la_1}{\la_0}\rme^{-2\la t}  &  \quad     1-\rme^{-2\la t}
\end{pmatrix},
\end{equation}
where $2\la=\la_0+\la_1.$ 

The density function of $X(t)$ is given by
\begin{equation}
\label{eq:PP}
\mathbf p^X(t,   \rmd y) 
=\rme^{-t\Lambda^{\mathrm{diag}}}\mathbf g*\mathbf q^t(\rmd y)
+\la_0\la_1\int_0^t\left[\mathcal B(t-\tau)\rme^{-\tau\Lambda^{\mathrm{diag}}}
\mathbf g*\mathbf q^\tau(\rmd y)\right]\rmd\tau,
\end{equation}
$t>0,$
where vector $\mathbf g*\mathbf q^t(\rmd y)$ is determined  by the entries
$g_i*q_i^t(\rmd y),\;i\in\{0, 1\};$ see \eqref{def:dnXt-renewal}.
\end{theo}

\proof Due to \eqref{def:dnXt-renewal}
by conditioning on the first switching similarly to \eqref{eq:eq-telegraph}
we have the coupled integral equations
\begin{equation*}
\left\{
\begin{aligned}
  p_0^X(t, \rmd y)=  & \rme^{-\la_0t}g_0*q_0^t(\rmd y)
  +\int_0^t\la_0\rme^{-\la_0\tau}p_1^X(t-\tau, \rmd y)\rmd \tau, \\
  p_1^X(t, \rmd y)=&  \rme^{-\la_1t}g_1*q_1^t(\rmd y)
  +\int_0^t\la_1\rme^{-\la_1\tau}p_0^X(t-\tau, \rmd y)\rmd\tau,
\end{aligned}
\right.\qquad t\geq0.
\end{equation*}

Therefore, the  time-Laplace transforms, 
\[\pi_i(s, \cdot):=\mathcal L_{t\to s}\left[p_i^X(t, \cdot)\right]=\int_0^\infty\rme^{-st}p_i^X(t, \cdot)\rmd t,\] 
satisfy the algebraic system,
\begin{equation}
\label{eq:P}
\left\{
\begin{aligned}
\pi_0(s, \cdot)=&Q_0(s, \cdot)+\frac{\la_0}{\la_0+s}\pi_1(s, \cdot),\\
\pi_1(s, \cdot)=&Q_1(s, \cdot)+\frac{\la_1}{\la_1+s}\pi_0(s, \cdot),
\end{aligned}
\right.
\end{equation}
where $Q_i(s, \cdot)$ is the time-Laplace transform of  $\rme^{-\la_it}g_i*q_i^t(\cdot):$
\[
Q_i(s, \cdot)=\mathcal L_{t\to s}\left[
\rme^{-\la_it}g_i*q_i^t(\cdot)
\right],\qquad i\in\{0, 1\}.
\]
The solution of \eqref{eq:P} can be written in the form
\begin{equation*}
\left\{
\begin{aligned}
\pi_0(s, \cdot)=  Q_0(s, \cdot)+&\frac{\la_0\la_1}{2\la}\left[Q_0(s, \cdot)+Q_1(s, \cdot)\right]\frac{1}{s}\\
   +&\frac{\la_0}{2\la}\left[\la_0Q_1(s, \cdot)-\la_1Q_0(s, \cdot)\right]\frac{1}{2\la+s},  \\
\pi_1(s, \cdot)=  Q_1(s, \cdot)+&\frac{\la_0\la_1}{2\la}\left[Q_0(s, \cdot)+Q_1(s, \cdot)\right]\frac{1}{s}\\
   -&\frac{\la_1}{2\la}\left[\la_0Q_1(s, \cdot)-\la_1Q_0(s, \cdot)\right]\frac{1}{2\la+s}.
\end{aligned}
\right.
\end{equation*}

By applying the inverse Laplace transformation $\mathcal L^{-1}_{s\to t}$
we obtain \eqref{eq:PP}, cf. \cite{STAPRO131}.

Theorem \ref{theo:renewal} allows to obtain the Fourier transform 
$\EE\left[\rme^{\mathrm i\theta X(t)}\right]$
of $X(t).$ 

Let 
$\psi_0=\psi_0(\theta)$ and $\psi_1=\psi_1(\theta)$ be the L\'evy-Khintchine exponents 
of the underlying L\'evy blocks $\eta_0$ and $\eta_1$ respectively
and
\[
\Psi^{\mathrm{diag}}(\theta)=
\begin{pmatrix}
     \psi_0(\theta) & \quad  0 \\ \\
   0   &  \quad \psi_1(\theta)
\end{pmatrix}.
\]
Note that the Fourier transform of the convolution $\mathbf g*\mathbf q^t(\rmd x)$ is given by
\[
\mathcal F_{x\to\theta}[\mathbf g*\mathbf q^t(\rmd x)]
=\int_{-\infty}^\infty\rme^{\mathrm i\theta x}\mathbf g*\mathbf q^t(\rmd x)
=\rme^{-t\Psi^{\mathrm{diag}}(\theta)}\widehat{\mathbf g}(\theta),
\] 
where $\widehat{\mathbf g}(\theta)=(\widehat g_0(\theta),\;\widehat g_1(\theta))^\ttT,$
\[
\widehat g_0(\theta)=\int_{-\infty}^\infty\rme^{\mathrm i\theta x}g_0(\rmd x),\qquad
\widehat g_1(\theta)=\int_{-\infty}^\infty\rme^{\mathrm i \theta x}g_1(\rmd x),
\]
are  the Fourier transforms of the distributions of starting points.

Therefore,
by \eqref{eq:PP}   the Fourier transform
$\vec{\Phi}=(\mathrm\Phi_0, \mathrm\Phi_1)^{\tt T}$ of 
$p^X(t, \cdot),$  
\[
\mathrm\Phi_i(t, \theta)=\EE\left(
\rme^{\mathrm i\theta X(t)}~|~\ep(0)=i
\right),\qquad i\in\{0, 1\},  \quad t>0,
\]
is expressed by
\begin{equation}\label{eq:phi-M}
\vec\Phi(t, \theta)=\left[
\rme^{-t(\Lambda^{\mathrm{diag}}+\Psi^{\mathrm{diag}}(\theta))}
+\la_0\la_1\int_0^t\mathcal B(t-\tau)
\rme^{-\tau(\Lambda^{\mathrm{diag}}+\Psi^{\mathrm{diag}}(\theta))}\rmd\tau
\right]\widehat{\mathbf g}(\theta).
\end{equation}

Integrating in  \eqref{eq:phi-M}  by  \eqref{def:B}
we obtain the following   formulae for the entries 
$\mathrm\Phi_0(t, \theta)$ and $\mathrm\Phi_1(t, \theta)$ of $\vec \Phi(t, \theta)$.
\begin{theo}\label{theo:charfun2}
The Fourier  transform $\vec \Phi(t, \theta)$ of $X(t)$  is given  by
\begin{equation}\label{eq0:charfun2}
\mathrm\Phi_0(t, \theta)=\rme^{-t(\la_0+\psi_0(\theta))}\widehat g_0(\theta)
+\frac{\la_0\la_1}{2\la}\left\{ \mathrm\Phi_{00}(t, \theta)\widehat g_0(\theta)
+\mathrm\Phi_{01}(t, \theta)\widehat g_1(\theta)
\right\}
\end{equation}
and
\begin{equation}\label{eq1:charfun2}
\mathrm\Phi_1(t, \theta)=\rme^{-t(\la_1+\psi_1(\theta))}\widehat g_1(\theta)
+\frac{\la_0\la_1}{2\la}\left\{ \mathrm\Phi_{10}(t, \theta)\widehat g_0(\theta)
+\mathrm\Phi_{11}(t, \theta)\widehat g_1(\theta)
\right\},
\end{equation}
where 
\[
\mathrm\Phi_{00}(t, \theta)=
\frac{1-\exp\left(-(\la_0+\psi_0(\theta))t\right)}{\la_0+\psi_0(\theta)}
-\frac{ \exp(-2\la t)-\exp\left(-(\la_0+\psi_0(\theta))t\right)}{\psi_0(\theta)-\la_1},
\]
\[
\mathrm\Phi_{01}(t, \theta)=
\frac{1-\exp\left(-(\la_1+\psi_1(\theta))t\right)}{\la_1+\psi_1(\theta)}
+\frac{\la_0}{\la_1}\frac{ \exp(-2\la t)-\exp\left(-(\la_1+\psi_1(\theta))t\right)}{\psi_1(\theta)-\la_0},
\]
\[
\mathrm\Phi_{10}(t, \theta)=
\frac{1-\exp\left(-(\la_0+\psi_0(\theta))t\right)}{\la_0+\psi_0(\theta)}
+\frac{\la_1}{\la_0}\frac{ \exp(-2\la t)-\exp\left(-(\la_0+\psi_0(\theta))t\right)}{\psi_0(\theta)-\la_1},
\]
and
\[
\mathrm\Phi_{11}(t, \theta)=
\frac{1-\exp\left(-(\la_1+\psi_1(\theta))t\right)}{\la_1+\psi_1(\theta)}
-\frac{ \exp(-2\la t)-\exp\left(-(\la_1+\psi_1(\theta))t\right)}{\psi_1(\theta)-\la_0}.
\]
\end{theo}

Formulae \eqref{eq0:charfun2}-\eqref{eq1:charfun2} allows us to prove the limit theorem.

To start with,  define a simple mixture of distributions. 
Let \[F^{X_1}(\rmd x), \ldots, F^{X_n}(\rmd x)\] 
be the distributions of random variables $X_1, \ldots, X_n,$ and $\zeta$
be an independent random variable taking values $1,\ldots, n$ with probabilities 
$p_1, \ldots, p_n,$ $\sum_{k=1}^np_k=1.$ 
The mixture of distributions, $\mathfrak M_\bp(X_1,\;\ldots,\; X_n),$ 
is defined by 
\[
\mathfrak M_\bp(X_1,\;\ldots,\; X_n)\stackrel{D}{=}X_\zeta,\qquad \bp=(p_1,\ldots,p_n).
\]
The distribution of $X_\zeta$ is given by linear combination,
$
\;\sum_{k=1}^np_kF^{X_k}(\rmd x).
$

For example,  an asymmetric  double exponential distribution
 (an asymmetric   Laplace distribution)
with the density 
\[
f(x)=p\la_0\rme^{-\la_0x}\1_{\{x>0\}}+(1-p)\la_1\rme^{-\la_1x}\1_{\{x<0\}},
\]
  see e. g. \cite{CaiKou}, is the mixture of $X_0$ and $-X_1$, where
  $X_0$  and $X_1$ are independent and 
  exponentially distributed, $X_i\sim\mathrm{Exp}(\la_i),\;i\in\{0, 1\}.$
 
\begin{theo}\label{theo:lim}
Let $X=X(t),\;t\geq0,$ be defined by \eqref{eq2:X} and 
$\tau_0, \;\tau_1$ are independent
exponentially distributed random variables with parameters $\la_0$ and $\la_1$ 
respectively\textup{.} 

$X(t)$ converges in distribution \textup{(}when $t\to\infty)$
to the $\bp$-mixture\textup{,} 
\begin{equation}
\label{theo:convergence}
X(t)\stackrel{D}{\longrightarrow}\mathfrak M_\bp(x_0+\eta_0(\tau_0),\;x_1+\eta_1(\tau_1)),
\end{equation}
where $\bp=(p_0,\;p_1),\;p_0=\dfrac{\la_1}{\la_0+\la_1},$ $p_1=\dfrac{\la_0}{\la_0+\la_1};\;$
$x_0,\; x_1$ are the independent
 renewal starting points with distributions $g_0,\;g_1.$ 

The Fourier  transform of the limiting distribution  
$\mathfrak M_\bp(x_0+\eta_0(\tau_0),\;x_1+\eta_1(\tau_1))$
is given by
\begin{equation}
\label{eq:limit_FT}
\mathrm\Phi_\infty(\theta):=\frac{\la_0\la_1}{\la_0+\la_1}
\left(\frac{\widehat g_0(\theta)}{\la_0+\psi_0(\theta)}
+\frac{\widehat g_1(\theta)}{\la_1+\psi_1(\theta)}\right).
\end{equation}
\end{theo}

\proof Note that $\mathrm{Re}[\psi_0(\theta)]\geq0,\;\mathrm{Re}[\psi_1(\theta)]\geq0$,
see \eqref{eq:Repsi}.
Passing  to limit  in \eqref{eq0:charfun2}-\eqref{eq1:charfun2} as $t\to\infty$
one can easily obtain $\lim_{t\to\infty}\mathrm\Phi_0(t, \theta)
=\lim_{t\to\infty}\mathrm\Phi_1(t, \theta)=\mathrm\Phi_\infty(\theta),$
where $\mathrm\Phi_\infty$ is given by \eqref{eq:limit_FT}.  Since,
\begin{equation*}
\begin{aligned}
\EE\left[\rme^{\mathrm i\theta\eta_i(\tau_i)}\right]
=&\int_0^\infty\la_i\rme^{-\la_i\tau}\EE\left[\rme^{\mathrm i\theta \eta_i(\tau)}\right]\rmd\tau\\
=&\int_0^\infty\la_i\rme^{-(\la_i+\psi_i(\theta))\tau}\rmd\tau
=\frac{\la_i}{\la_i+\psi_i(\theta)},\qquad
i\in\{0, 1\},
\end{aligned}\end{equation*}
the Fourier transform 
$\mathrm\Phi_\infty$ corresponds to the mixture \[\mathfrak M_\bp(x_0+\eta_0(\tau_0),\;x_1+\eta_1(\tau_1)).\]

\begin{ex}\textbf{Telegraph process.}
Let $\eta_0(t)=c_0t,\;\eta_1(t)=c_1t$ be deterministic\textup{.} 
The limiting distribution of the Goldstein-Kac process $X(t)$ 
with renewal starting points
is determined by \eqref{eq:limit_FT}
with $\psi_i(\theta)=-\mathrm ic_i\theta,\;i\in\{0, 1\}.$

Note that the inverse Fourier transform of 
$\dfrac{1}{\alpha-\mathrm ic\theta},\;\alpha>0,$ is given by
\begin{equation}\label{inv2}
\mathcal F^{-1}_{\theta\to x}\left[\frac{1}{\alpha-\mathrm ic\theta}\right]=
\left\{
\begin{aligned}
\frac{1}{|c|}\rme^{-\alpha x/c}\1_{\{cx>0\}},&\qquad c\neq0,   \\
    \frac{1}{\alpha}\delta(x),&  \qquad c=0.
\end{aligned}
\right.
\end{equation}

In the case of a particle always restarting from the origin\textup{,} that is
$g_0=g_1=\delta(\rmd x),$  the limit \eqref{theo:convergence}
is the mixture of two exponential distributions. 
In particular\textup{,}
if $c_0, c_1>0,$ 
then by  \eqref{inv2}  the limit   is characterised by
\[
f^*(x)=\frac{\la_0\la_1}{\la_0+\la_1}\left[
\frac{1}{c_0}\rme^{-\la_0x/c_0}+\frac{1}{c_1}\rme^{-\la_1x/c_1}
\right]\1_{\{x>0\}}.
\]
 If $c_0, c_1<0,$ 
then the limiting density function is 
\[
f^*(x)=\frac{\la_0\la_1}{\la_0+\la_1}\left[
\frac{1}{|c_0|}\rme^{-\la_0x/c_0}+\frac{1}{|c_1|}\rme^{-\la_1x/c_1}
\right]\1_{\{x<0\}}.
\]
In the case of the opposite signs\textup{,} $c_0>0>c_1,$ the limit  is 
the Laplace distribution with the density function
\[
f^*(x)=\frac{\la_0\la_1}{\la_0+\la_1}\left[
\frac{1}{c_0}\rme^{-\la_0x/c_0}\1_{\{x>0\}}+\frac{1}{|c_1|}\rme^{-\la_1x/c_1}
\1_{\{x<0\}}\right],
\]
cf \cite{STAPRO131}. If $c_0=0$ or $c_1=0,$ 
then the corresponding term becomes Dirac\textup{'}s $\delta(x)$.
\end{ex}

\begin{ex}\textbf{Alternating Brownian motions with  alternating drifts.}\label{ex:Brown}
Let $\eta_0(t)=c_0t+\sigma_0W^{(0)}(t),$ $\eta_1(t)=c_1t+\sigma_1W^{(1)}(t),$ where 
$W^{(0)}$ and $W^{(1)}$ are two independent standard Brownian motions\textup{,} 
$c_0, c_1, \sigma_0, \sigma_1$ are constants. In this case 
$\psi_0(\theta)=-\mathrm ic_0\theta+\sigma_0^2\theta^2/2,\;$ 
$\psi_1(\theta)=-\mathrm ic_1\theta+\sigma_1^2\theta^2/2$ and \eqref{eq:limit_FT} becomes 
\[
\frac{\la_0\la_1}{\la_0+\la_1}\left[
\frac{\widehat g_0(\theta)}{\sigma_0^2\theta^2/2-\mathrm ic_0\theta+\la_0}
+\frac{\widehat g_1(\theta)}{\sigma_1^2\theta^2/2-\mathrm ic_1\theta+\la_1}
\right].
\]
Note that
\[
\frac{1}{\sigma^2\theta^2/2-\mathrm ic\theta+\la}
=A\left(\frac{1}{\alpha_2-\mathrm i\theta}-\frac{1}{\alpha_1-\mathrm i\theta}\right),
\]
where $A>0$ and the numbers $\alpha_1, \alpha_2$ 
are of the opposite signs\textup{:}
\begin{equation}\label{def:Aalpha}
A=(c^2+2\la\sigma^2)^{-1/2},\quad 
\alpha_1=\frac{-c-\sqrt{c^2+2\la\sigma^2}}{\sigma^2}<0<
\frac{-c+\sqrt{c^2+2\l
a\sigma^2}}{\sigma^2}=\alpha_2.
\end{equation}
Hence\textup{,} by    \eqref{inv2}  
the inverse Laplace transform is
\[
\mathcal F^{-1}_{\theta\to x}\left[\frac{1}{\sigma^2\theta^2/2-\mathrm ic\theta+\la}\right]
=A\left[\rme^{-\alpha_2x}\1_{\{x>0\}}+\rme^{-\alpha_1x}\1_{\{x<0\}}
\right],
\]
which gives the Laplace distribution.

\begin{figure}[h]
\begin{center}
  \includegraphics[scale=0.66]{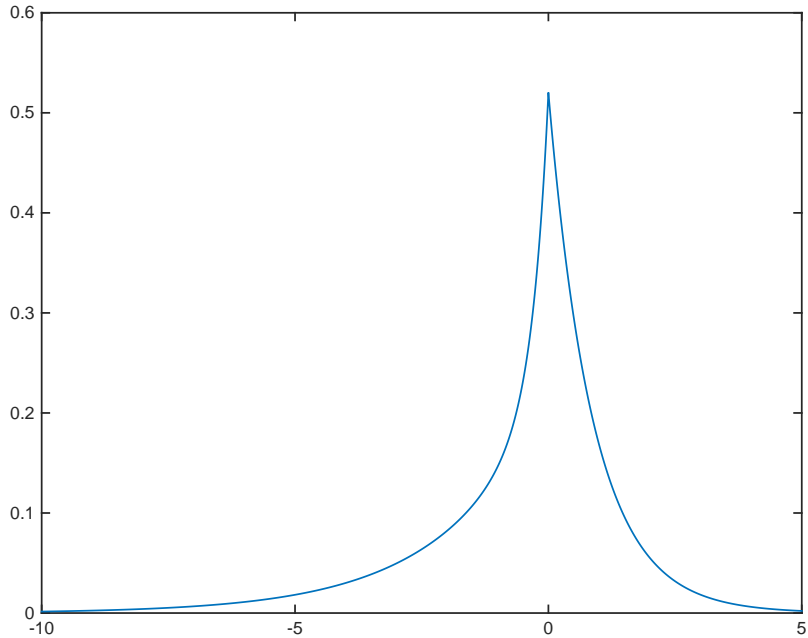}\\
 \caption{Density function $f^*(x),$ \eqref{eq:ex2}, of Example \ref{ex:Brown}
 with $\lambda_0=2,\; \la_1=1;\;c_0=1,\; c_1=-1;\;\sigma_0=1,\; \sigma_1=2$.
 {}{}}\label{fig1}
 \end{center}
\end{figure} 

Therefore\textup{,} the limiting distribution is given by the density function
\begin{equation*}
\begin{aligned}
f^*(x)=\frac{\la_0\la_1}{\la_0+\la_1}&\bigg\{
A^{(0)}\left[
\int_{-\infty}^x\rme^{-\alpha_2^{(0)}(x-y)}g_0(\rmd y)
+\int_{x}^\infty\rme^{-\alpha_1^{(0)}(x-y)}g_0(\rmd y)
\right]\\
&+A^{(1)}\left[
\int_{-\infty}^x\rme^{-\alpha_2^{(1)}(x-y)}g_1(\rmd y)
+\int_{x}^\infty\rme^{-\alpha_1^{(1)}(x-y)}g_1(\rmd y)
\right]
\bigg\},
\end{aligned}
\end{equation*}
where $A^{(0)},\;\alpha_1^{(0)}, \alpha_2^{(0)}$ and $ A^{(1)},\;\alpha_1^{(1)}, \alpha_2^{(1)}$
are defined by \eqref{def:Aalpha} with the sets of parameters $c_0, \sigma_0, \la_0$ and 
$c_1, \sigma_1, \la_1$ respectively.

If $g_0=g_1=\delta_0(\rmd x),$ then the limiting density function is given by 
the mixture of two Laplace distributions
\begin{align}
f^*(x)=\frac{\la_0\la_1}{\la_0+\la_1}\Big[\Big(A^{(0)}\rme^{-\alpha_2^{(0)}x}+&
A^{(1)}\rme^{-\alpha_2^{(1)}x}\Big)\1_{\{x>0\}}\nonumber\\
&+\left(A^{(0)}\rme^{-\alpha_1^{(0)}x}+
A^{(1)}\rme^{-\alpha_1^{(1)}x}\right)\1_{\{x<0\}}\Big],\label{eq:ex2}
\end{align}
see Fig. \ref{fig1}.
\end{ex}

\begin{ex}\textbf{Alternating compound Poisson processes with  
exponentially distributed jumps.}
Let $\eta_0=c_0t+\sum_{m=1}^{M_0(t)}Y_m^{(0)}$ 
and  $\eta_1=c_1t+\sum_{m=1}^{M_1(t)}Y_m^{(1)},$ 
where $\{Y_m^{(0)}\}, \;  \{Y_m^{(1)}\}$ are 
two independent sets of iid random variables\textup{,} 
independent of two independent counting Poisson processes 
$M_0(t),\;M_1(t)$ 
with intensities $\nu_0$ and $\nu_1$ respectively.
Let $Y_m^{(0)}\sim \mathrm{Exp}(a_0)$ and 
$Y_m^{(1)}\sim \mathrm{Exp}(a_1),\;a_0, a_1>0.$

The L\'evy-Khintchine exponents are
\[
\psi_i(\theta)=-\mathrm i\theta\left(c_i+\frac{\nu_i}{a_i-\mathrm i\theta}\right),\qquad i\in\{0, 1\},
\]
and \eqref{eq:limit_FT} becomes 
\[\begin{aligned}
\frac{\la_0\la_1}{\la_0+\la_1}\Big[
\frac{a_0-\mathrm i\theta}{a_0\la_0-\mathrm i\theta(a_0c_0+\nu_0+\la_0)-c_0\theta^2}
\widehat g_0(\theta)&\\
+
\frac{a_1-\mathrm i\theta}{a_1\la_1-\mathrm i\theta(a_1c_1+\nu_1+\la_1)-c_1\theta^2}&
\widehat g_1(\theta)
\Big].
\end{aligned}\]

In this case, we manage like in Example \ref{ex:Brown}.
Note that
\[
\frac{a-\mathrm i\theta}{a\la-\mathrm i\theta(ac+\nu+\la)-c\theta^2}
=\left\{\begin{aligned}
\frac{A_1}{\alpha_1-\mathrm i\theta}+\frac{A_2}{\alpha_2-\mathrm i\theta},&
\qquad\text{if~~~}c\neq0,\\
\frac{1}{\nu+\la}+\frac{A_3}{\alpha_3-\mathrm i\theta},&\qquad\text{if~~~}c=0.
\end{aligned}\right.
\]
Here $\alpha_{1}$ and  $\alpha_{2}$ are the roots of the
equation $\phi(\alpha):=c\alpha^2-(ac+\nu+\la)\alpha+a\la=0,\;c\neq0.$
Since $\phi(0)=a\la>0$ and $\phi(a)=-\nu a<0$ we have two positive roots $\alpha_1,\;\alpha_2>0,$
if $c>0;$ in the case $c<0$ the roots are of the opposite signs, $\alpha_2<0<\alpha_1<a$.
Precisely,
\begin{equation}\label{eq:alphak}
\begin{aligned}
\alpha_1=\frac{ac+\nu+\la-\sqrt{D}}{2c},&\qquad 
\alpha_2=\frac{ac+\nu+\la+\sqrt{D}}{2c},\\
\alpha_3=&a\la/(\nu+\la)>0
\end{aligned}\end{equation}
 and
\begin{equation}\label{eq:Ak}
\begin{aligned}
A_1=\frac{a-\alpha_1}{\sqrt{D}},&\qquad A_2=\frac{\alpha_2-a}{\sqrt{D}},\\
A_3=(a-\alpha_3)/(\nu+\la)=&a\nu/(\nu+\la)^2>0,
\end{aligned}\end{equation}
where $D=(ac+\nu+\la)^2-4ac\la,\;D>0.$
Note that if $c>0,$  then $A_1,\;A_2>0;$  if $c<0,$  then $A_1>0>A_2$. 
%

Due to \eqref{inv2} the inverse Fourier transform 
\[F(x; a, c, \nu, \la)=\mathcal F^{-1}_{\theta\to x}
\left[\frac{a-\mathrm i\theta}{a\la-\mathrm i\theta(ac+\nu+\la)-c\theta^2}\right]\]
can be expressed explicitly, 
\[
F(x; a, c, \nu, \la)=\left\{\begin{aligned}
\left[A_1\rme^{-\alpha_1x}+A_2\rme^{-\alpha_2x} 
\right]\1_{\{x>0\}},&\qquad\text{if~~~}c>0,\\
A_1\rme^{-\alpha_1x}\1_{\{x>0\}}-A_2\rme^{-\alpha_2x}\1_{\{x<0\}},&\qquad \text{if~~~}c<0,\\
\frac{\delta(x)}{\nu+\la}+A_3\rme^{-\alpha_3x}\1_{\{x>0\}},&\qquad\text{if~~~}c=0,
\end{aligned}\right.
\]
where $\alpha_k$ and $A_k,\;k=1, 2, 3,$ are defined by
\eqref{eq:alphak} and \eqref{eq:Ak}.

So the limiting density function is given by the convolutions
\[f^*(x)=\frac{\la_0\la_1}{\la_0+\la_1}
\left[F(\cdot;  a_0, c_0, \nu_0, \la_0)*g_0(x)+F(\cdot;  a_1, c_1, \nu_1, \la_1)*g_1(x)\right].\]

The case of $g_0=g_1=\delta_0(\rmd x)$ can be described 
as in the previous  examples.
\end{ex}

\section{Kac-L\'evy processes with jumps}\label{sec:jump}
\setcounter{equation}{0}

Unlike the case of Section \ref{sec:renewal}
suppose that the Markov-modulated 
process $X=X(t)$ (with alternating intensities $\la_0, \la_1$) is jumping after each switching:
\begin{align}
X(t)=&X(T_n-)+Y_{n-1}+\eta_n(t-T_n),\qquad &T_n\leq t<T_{n+1},&\qquad n\geq1,\nonumber\\
X(t)=&\eta_0(t),\qquad &0\leq t<T_1,&\label{def:X-jump}
\end{align}
where the underlying L\'evy processes $\{\eta_n(t),\;t\geq0\}_{n\geq0}$ 
and 
jump magnitudes $\{Y_n\}_{n\geq0}$ are independent with  alternating distributions.
Such Markov-modulated L\'evy process $X$ can be considered as an example of 
It\^o-L\'evy process, see e.g. \cite{DiNunnoOksendal}.

\subsection{L\'evy-Laplace exponent of $X(t)$}
Let
$\ell_0$ and $\ell_1$ be alternating L\'evy-Laplace exponents of the underlying processes $\eta_n(t)$
and $h_0(\rmd y)$ and $h_1(\rmd y)$ be alternating distributions of jump magnitudes $Y_n$ with 
the Laplace transforms 
\[
\widetilde h_0(\xi)=\int_{-\infty}^\infty\rme^{-\xi y}h_0(\rmd y),\qquad  
\widetilde h_1(\xi)=\int_{-\infty}^\infty\rme^{-\xi y}h_1(\rmd y).
\]

We denote  
\begin{align}\label{def:lamu}
\la:=(\la_0+\la_1)/2,\quad& \mu:=(\la_0-\la_1)/2,\\
\label{def:ellemm}
 \ell(\xi):=(\ell_0+\ell_1)/2,\quad& \emm(\xi):=(\ell_0-\ell_1)/2.
\end{align}

Let 
\begin{equation}\label{matrixA:a}
\mathcal L(\xi)=\begin{pmatrix}
      \la_0+\ell_0(\xi)& \quad-\la_0\widetilde h_0(\xi)   \\ \\
     -\la_1\widetilde h_1(\xi) &  \quad      \la_1+\ell_1(\xi)
\end{pmatrix}.
\end{equation}
Similar to \eqref{eq:eq-telegraph}, 
the Laplace transform $\vec\Ell(t, \xi)=(\Ell_0(t, \xi),\;\Ell_1(t, \xi))^\ttT$ of $X(t),$ 
\[
\Ell_i(t, \xi)=\EE\left[\rme^{-\xi X(t)}~|~\ep(0)=i\right],\qquad i\in\{0, 1\},
\]
can be obtained explicitly.

\begin{theo}\label{theo:charfun1}
For $t\geq0$ 
\begin{align}
\label{eq0:charfun1}
    \Ell_0(t, \xi)=\frac12\rme^{-(\la+\ell(\xi))t}\bigg[
    \rme^{tD(\xi)}+&\rme^{-tD(\xi)}\\
    +&
    \frac{\la_0\widetilde h_0(\xi)-\mu-\emm(\xi)}{D(\xi)}\left(
     \rme^{tD(\xi)}-\rme^{-tD(\xi)}
    \right)\nonumber
    \bigg],   \\ \nonumber\\
    \label{eq1:charfun1}
       \Ell_1(\xi, t)=\frac12\rme^{-(\la+\ell(\xi))t}\bigg[
        \rme^{tD(\xi)}+&\rme^{-tD(\xi)}\\+&
        \frac{\la_1\widetilde h_1(\xi)+\mu+\emm(\xi)}{D(\xi)}\left(
        \rme^{tD(\xi)}-\rme^{-tD(\xi)}
        \right)\nonumber
    \bigg],
\end{align}
where 
$D(\xi)=\left[(\emm(\xi)+\mu)^2+\la_0\la_1\widetilde h_0(\xi)\widetilde h_1(\xi)\right]^{1/2}>0.$ 

Here $\la, \mu$ and $\ell(\xi), \emm(\xi)$
are defined by \eqref{def:lamu} and \eqref{def:ellemm}.
\end{theo}

\proof
A commonplace is that
\[ 
\vec\Ell(t, \xi)=\exp\left(-t\mathcal L(\xi)\right)\mathbf 1, \qquad \mathbf 1=(1, 1)^\ttT,
\] 
where $\mathcal L(\xi)$ is defined by \eqref{matrixA:a},
cf. \cite[Chap. XI, Proposition 2.2]{asmussen}.
Therefore,
\begin{equation}\label{eq2:Lexp}
\vec\Ell(t, \xi)=\exp\left(-t\mathcal L(\xi)\right)\mathbf 1
=\rme^{-\alpha_1t}\mathbf{e_1}+\rme^{-\alpha_2t}\mathbf{e_2},\qquad t\geq0,
\end{equation}
 where 
$\alpha_1$ and $\alpha_2$ are the eigenvalues of the matrix $\mathcal L(\xi)$
and $\mathbf {e_1},\;\mathbf {e_2}$ are corresponding eigenvectors 
with $\mathbf {e_1}+\mathbf {e_2}=\mathbf 1.$

Representation \eqref{eq2:Lexp}  can be derived as a consequence of the  
coupled integral equations arising by conditioning on the first pattern's switch,
 \begin{equation}\label{eq:p0Xp1X}
 \left\{
\begin{aligned}
   p_0^X(t, \cdot) &= \rme^{-\la_0t}q_0^t(\cdot)+\int_0^t\la_0\rme^{-\la_0\tau}
   \left[p_1^X(t-\tau, \cdot)*q_0^\tau* h_0\right](\cdot)\rmd\tau,  \\
    p_1^X(t, \cdot) &= \rme^{-\la_1t}q_1^t(\cdot)+\int_0^t\la_1\rme^{-\la_1\tau}
   \left[p_0^X(t-\tau, \cdot)*q_1^\tau* h_1\right](\cdot)\rmd\tau,
\end{aligned}\right.
\end{equation}
where the convolution $*$ is defined by
\[
\int_{\RR_z}\varphi(z) \left[ p(t-\tau, \cdot)*q^\tau* h\right](\rmd z)
=\int_{\RR^3}\varphi(x+y+z) p(t-\tau, \rmd x)q^\tau(\rmd y)h(\rmd z)
\]
for any test-function $\varphi$. 
From \eqref{eq:p0Xp1X}
one can obtain  
the following coupled integral equations for the Laplace transform of $X(t),$
\begin{equation}\label{eq1:phi0phi1}
\left\{
\begin{aligned}
    \Ell_0(t, \xi)=&\exp\left\{-(\la_0+\ell_0(\xi))t\right\}
+\la_0\widetilde h_0(\xi)
\int_0^t\rme^{-\la_0\tau}\rme^{-\ell_0(\xi)\tau}\Ell_1(t-\tau, \xi)\rmd \tau,   \\
     \Ell_1(t, \xi)=&\exp\left\{-(\la_1+\ell_1(\xi))t\right\}
+\la_1\widetilde h_1(\xi)
\int_0^t\rme^{-\la_1\tau}\rme^{-\ell_1(\xi)\tau}\Ell_0(t-\tau, \xi)\rmd \tau, 
\end{aligned}\right.
\end{equation}
which is equivalent to 
the   initial value problem:
\begin{equation}
\label{eq2:phi0phi1}
\frac{\rmd\vec\Ell}{\rmd t}(t, \xi)=-\mathcal L(\xi)\vec\Ell(t, \xi),\qquad t>0,
\end{equation}
\[
\vec\Ell(0, \xi)=\mathbf 1,
\]
where 
$\mathcal L(\xi)$ is defined by \eqref{matrixA:a}. 
Formula \eqref{eq2:Lexp}  is the solution of  \eqref{eq2:phi0phi1}.

The eigenvalues $\alpha_1,\;\alpha_2$ and the corresponding eigenvectors
$\mathbf {e_1},\;\mathbf {e_2}$ of the matrix $\mathcal L(\xi)$ can be found explicitly. 
We have 
\begin{equation}\label{def:alpha12}
\alpha_{1}=\lambda+\ell(\xi)- D(\xi),\qquad \alpha_{2}=\lambda+\ell(\xi)+ D(\xi)
\end{equation}
and  
\begin{equation}\label{def:e1}
\mathbf {e_1}=\frac12
\left(
1-\frac{\mu+\emm(\xi)-\la_0\widetilde h_0(\xi)}{D(\xi)},\; 
1+\frac{\mu+\emm(\xi)+\la_1\widetilde h_1(\xi)}{D(\xi)}
\right)^{\tt T},
\end{equation} 
\begin{equation}\label{def:e2}
\mathbf {e_2}=\frac12
\left(
1+\frac{\mu+\emm(\xi)-\la_0\widetilde h_0(\xi)}{D(\xi)},\; 
1-\frac{\mu+\emm(\xi)+\la_1\widetilde h_1(\xi)}{D(\xi)}
\right)^{\tt T}.
\end{equation}

 Formulae \eqref{eq0:charfun1}-\eqref{eq1:charfun1} follow 
 from \eqref{eq2:Lexp} and \eqref{def:alpha12}-\eqref{def:e2}.
 
 \begin{rem}
 Notice that in the case of Kac process which corresponds to the L\'evy-Laplace exponents
  $\ell_0=c_0\xi,\; \ell_1=c_1\xi$ formulae \eqref{eq0:charfun1}-\eqref{eq1:charfun1} of
 Theorem \ref{theo:charfun1} are equivalent to \cite[(3.3), Theorem 3.1]{STAPRO82},
 where the characteristic functions of $X(t)$ are presented explicitly. 
 \end{rem}

Assume that
the Markov modulation $\ep=\ep(t)$ is determined by the underlying L\'evy processes. 
For example, let $\ep(t)$ be
switching at times of \textup{``}big\textup{"} jumps.   
Let  $\ell_0$ switches to $\ell_1$ just after 
a \textup{``}big\textup{"} negative jump of the current L\'evy pattern $\eta$, 
that is if $\Delta\eta_0(t)<-R_0;$ vice versa,
$\ell_1$ switches to $\ell_0$ just after 
a \textup{``}big\textup{"} positive jump\textup{,} that is if $\Delta\eta_1(t)>R_1$.

\begin{cor}
The Laplace transform $\vec\Ell=(\Ell_0, \Ell_1)^{\tt T}$ 
 of this process  can be obtained in the following form\textup{:}
\begin{align}
\label{eq:phi0}
    \Ell_0(t, \xi)=&\frac12\rme^{-(a(\xi)+\ell(\xi))t} 
    \left[\rme^{tD(\xi)}+\rme^{-tD(\xi)}
    +\frac{\la_0a_0(\xi)-b(\xi)- \emm(\xi)}{D(\xi)}\left(
    \rme^{tD(\xi)}-\rme^{-tD(\xi)}
    \right)
    \right],   \\
    \label{eq:phi1}
    \Ell_1(t, \xi)=&\frac12\rme^{-(a(\xi)+\ell(\xi))t}   
     \left[\rme^{tD(\xi)}+\rme^{-tD(\xi)}
    +\frac{\la_1a_1(\xi)+b(\xi)+\emm(\xi)}{D(\xi)}\left(
    \rme^{tD(\xi)}-\rme^{-tD(\xi)}
    \right)
    \right], \end{align}
Here $\ell(\xi)$ and $\emm(\xi)$ are defined by \eqref{def:ellemm} and
$D(\xi)^2=(\emm(\xi)+b(\xi))^2+\la_0\la_1a_0(\xi)a_1(\xi),$ where 
\begin{equation*}
\la_0=\mathrm\Pi_0\{x:~x<-R_0\},\qquad \la_1=\mathrm\Pi_1\{x:~x>R_1\},
\end{equation*}
\begin{equation*}
a_0(\xi)=\int_{x\leq -R_0}\rme^{-\xi x}\mathrm\Pi_0(\rmd x),
\qquad
a_1(\xi)=\int_{x\geq R_1}\rme^{-\xi x}\mathrm\Pi_1(\rmd x);
\end{equation*}
\begin{equation*}
a(\xi):=\frac12\left(a_0(\xi)+a_1(\xi)\right),\qquad
b(\xi):=\frac12\left(a_0(\xi)-a_1(\xi)\right).
\end{equation*}
\end{cor}

\proof
To prove it, note that in this case, instead of \eqref{eq1:phi0phi1}, we have equations
\begin{equation*}
\Ell_0(t, \xi)=\exp\{-\left(\la_0+\overline\ell_0(\xi)\right)t\}
 +\la_0a_0(\xi)
 \int_0^t\rme^{-\la_0\tau}\rme^{-\overline\ell_0(\xi)\tau}\Ell_1(t-\tau, \xi)\rmd \tau, 
\end{equation*}
\begin{equation*}
\Ell_1(t, \xi)=\exp\{-\left(\la_1+\overline\ell_1(\xi)\right)t\}
  +\la_1a_1(\xi)
  \int_0^t\rme^{-\la_1\tau}\rme^{-\overline\ell_1(\xi)\tau}\Ell_0(t-\tau, \xi)\rmd \tau.
\end{equation*}
Here $\la_0=\mathrm\Pi_0\{x:~x<-R_0\},\;\la_1=\mathrm\Pi_1\{x:~x>R_1\}$ and 
$\overline\ell_0(\xi),\;\overline\ell_1(\xi)$ are defined by
\begin{equation*}
\begin{aligned}
\overline\ell_0(\xi)=&c_0\xi-\frac12\sigma_0^2\xi^2
+\int_{-R_0}^\infty\left(1-\rme^{-\xi x}-\xi x\1_{\{x<R_0\}}
\right)\mathrm\Pi_0(\rmd x),\\
\overline\ell_1(\xi)=& c_1\xi-\frac12\sigma_1^2\xi^2
+\int_{-\infty}^{R_1}\left(1-\rme^{\xi x}-\xi x\1_{\{x>-R_1\}}
\right)\mathrm\Pi_1(\rmd x). 
\end{aligned}\end{equation*}
Further,
\[
\la_0+\overline\ell_0(\xi)\equiv a_0(\xi)+\ell_0(\xi),\qquad
\la_1+\overline\ell_1(\xi)\equiv a_1(\xi)+\ell_1(\xi).
\]
Representation \eqref{eq:phi0}-\eqref{eq:phi1} follows 
from \eqref{eq0:charfun1}-\eqref{eq1:charfun1}.

\subsection{Kac-L\'evy subordinated Kac-L\'evy process with jumps}
Let $S_0=S_0(t)$ and $S_1=S_1(t)$ be two subordinators (increasing L\'evy processes)
with the L\'evy-Laplace exponents 
$\ell_0^S(\xi)$ and $\ell_1^S(\xi)$ respectively.
Let $\ep^S=\ep^S(t)$ be a two-state Markov process with switching intensities $\la_0^S$ and $\la_1^S$.
We define a Markov-modulated subordinator as
 a Markov-modulated L\'evy process $Z$ based on $S_0,\;S_1$ and $\ep^S$.

Let $X=X(t)$ be a Kac-L\'evy process with jumps
based on 
the two-state Markov process $\ep^X$  
with switching intensities $\la_0^X$ and $\la_1^X$ and 
on the independent L\'evy processes $\eta_0$ and $\eta_1$
with the  L\'evy-Laplace exponents 
$\ell_0^X, \ell_1^X$. Assume that $X$ and  $Z$ are independent.

Consider the Kac-L\'evy subordinated process $X\circ Z(t).$ 
The  Laplace transform of   $X\circ Z(t),$ 
\[
\Ell_{ij}(t, \xi)=\EE\{\rme^{-\xi X\circ Z(t)}~|~\ep^X(0)=i,\;\ep^S(0)=j\},\qquad i, j\in\{0, 1\},
\]
can be expressed explicitly.

\begin{theo}
Let 
$
\vec\Ell_j(t, \xi)
=\left(\Ell_{0j}(t, \xi),\;\Ell_{1j}(t, \xi)\right),\;j\in\{0, 1\}.
$ Then
\begin{align}
   \vec L_0(t, \xi)=\rme^{-\left(\la^S+\ell^S(\alpha_1(\xi))\right)t}\left[
   \cosh(tD^S(\alpha_1(\xi)))+\frac{\la^S-\emm^S(\alpha_1(\xi))}{D^S(\alpha_1(\xi))}   \sinh(tD^S(\alpha_1(\xi)))
   \right]&\mathbf{e_1}\nonumber\\
+\rme^{-\left(\la^S+\ell^S(\alpha_2(\xi))\right)t}\left[
   \cosh(tD^S(\alpha_2(\xi)))+\frac{\la^S-\emm^S(\alpha_2(\xi))}{D^S(\alpha_2(\xi))}   \sinh(tD^S(\alpha_2(\xi)))
   \right]&\mathbf{e_2},\label{eq:Ell0}\\
    \vec L_1(t, \xi)=\rme^{-\left(\la^S+\ell^S(\alpha_1(\xi))\right)t}\left[
   \cosh(tD^S(\alpha_1(\xi)))+\frac{\la^S+\emm^S(\alpha_1(\xi))}{D^S(\alpha_1(\xi))}   \sinh(tD^S(\alpha_1(\xi)))
   \right]&\mathbf{e_1}\nonumber\\
+\rme^{-\left(\la^S+\ell^S(\alpha_2(\xi))\right)t}\left[
   \cosh(tD^S(\alpha_2(\xi)))+\frac{\la^S+\emm^S(\alpha_2(\xi))}{D^S(\alpha_2(\xi))}   \sinh(tD^S(\alpha_2(\xi)))
   \right]&\mathbf{e_2}.\label{eq:Ell1}
\end{align}
Here $D^S(\cdot)=\left[(\emm^S(\cdot)+\mu^S)^2+\la_0^S\la_1^S\right]^{1/2};$
$\alpha_i(\xi),\;\mathbf {e}_i,\;i\in\{0,1\}$ are identified by \eqref{def:alpha12}-\eqref{def:e2}.
\end{theo}
\proof
By  \eqref{eq2:Lexp}
\[
\EE\left(\rme^{-\xi X\circ Z(t)}~|~Z(t)\right)
=\rme^{-\alpha_1(\xi)Z(t)}\mathbf{e_1}(\xi)+\rme^{-\alpha_2(\xi)Z(t)}\mathbf{e_2}(\xi),
\]
where $\alpha_1, \alpha_2,\;\mathbf{e_1}, \mathbf{e_2}$ are defined by \eqref{def:alpha12}-\eqref{def:e2}.
Formulae \eqref{eq:Ell0}-\eqref{eq:Ell1} follow from \eqref{eq0:charfun1}-\eqref{eq1:charfun1}.

For the L\'evy subordination of L\'evy processes see  \cite[Theorem 30.1]{sato}. 

\subsection{Exponential functional }\label{sec:exponential}
Let $X$ be an alternating Markov-modulated 
L\'evy process with jumps based on a Markov process $\ep=\ep(t)\in\{0, 1\},$ 
with alternating switching intensities $\la_0,\;\la_1,$
the sequence of jumps $\{Y_n\}_{n\geq1}$ with alternating distributions
and  the independent alternating L\'evy processes $\{\eta_n\}_{n\geq0},$ \eqref{def:X-jump}.

We study the exponential functional
\begin{align}
I_\infty(\gamma)=&\int_0^{\infty}\rme^{-\gamma X(t)}\rmd t
=\sum_{n=0}^\infty\int_{T_n}^{T_{n+1}}\rme^{-\gamma X(t)}\rmd t\nonumber\\
=&\int_0^{T_1}\rme^{-\gamma\eta_0(t)}\rmd t+\sum_{n=1}^\infty\rme^{-\gamma X(T_n-)-\gamma Y_{n-1}}
\int_{T_n}^{T_{n+1}}\rme^{-\gamma\eta_n(t-T_n)}\rmd t, \label{def:Igamma}\\
=&\mathcal T_{0, \gamma}
+\sum_{n=1}^\infty(Z_n)^\gamma\cdot \mathcal T_{n, \gamma},\qquad \gamma>0,\nonumber
\end{align}
where
\[
\mathcal T_{n, \gamma}=\int_{T_n}^{T_{n+1}}\exp\left(-\gamma\eta_n(t-T_n)\right)\rmd t
=\int_0^{\Delta T_n}\exp\left(-\gamma \eta_n(t)\right)\rmd t,\]
$
\Delta T_n=T_{n+1}-T_n,
$
and $Z_n=\exp\left(-X(T_n-)-Y_{n-1}\right),\;n\geq0.$

By definition \eqref{def:X-jump}
\[
Z_n=\exp(-X(T_n-)-Y_{n-1})=\exp\left(-X(T_n)\right)
\]
\[
=\exp\left(-\sum_{k=1}^n\left[X(T_k)-X(T_{k-1})\right]\right)
=\exp\left(-\sum_{k=1}^{n}\left[\eta_{k-1}(\Delta T_{k})+Y_{k-1}\right]\right),\; n\geq1.
\]
Note that
$\mathcal T_{n, \gamma}$ are mutually independent and independent of $Z_n.$

Similar to \eqref{def:Igamma} the exponential functional
\begin{equation}
\label{def:I}
I_\infty=I_\infty(1)=\int_0^{\infty}\rme^{- X(t)}\rmd t.
\end{equation}
can be represented by
\[
I_\infty=\mathcal T_0+\sum_{n=1}^\infty Z_n\cdot\mathcal T_n,
\]
where
$
\mathcal T_n=\int_0^{\Delta T_n}\rme^{-\eta_n(t)}\rmd t,
\quad n\geq0.
$

First, we study the convergence of the exponential functional $I_\infty$. 
To begin, we need the following two important statements.

\begin{lem}\label{lem2}
Let $I_\infty(\gamma)$ be defined by \eqref{def:Igamma} 
with alternating L\'evy-Laplace exponents $\ell_0(\xi),\;\ell_1(\xi).$
If $\EE\left[I_\infty(\gamma)\right]<\infty$ for some $\gamma,\;\gamma\in(0, 1],$
such that
\begin{equation}
\label{eq:l01gamma01}
\ell_0(\gamma)+\la_0>0,\qquad \ell_1(\gamma)+\la_1>0,
\end{equation}
then
$I_\infty<\infty\;a.s.$ 
\end{lem}

\proof
 Let $\ep_n=\ep(T_n),\;n\geq0.$
If $\gamma$ satisfies \eqref{eq:l01gamma01}, then 
the expectations 
\[
\EE[\mathcal T_{n, \gamma}]
=\EE\left[\int_0^{\Delta T_n}\rme^{-\gamma\eta_n(t)}\rmd t\right]
=\int_0^\infty\la_{\ep_n}\rme^{-\la_{\ep_n}s}\left[\int_0^s\rme^{-t\ell_{\ep_n}(\gamma)}\rmd t\right]
\]
\[
=\frac{1}{\la_{\ep_n}+\ell_{\ep_n}(\gamma)}
\]
exist, positive and alternating. Let $\EE\left[I_\infty(\gamma)\right]<\infty.$ 
By \eqref{def:Igamma}
\[
\EE\left[I_\infty(\gamma)\right]=\EE\left[\mathcal T_{0, \gamma}\right]
+\sum_{n\geq1}\EE\left[(Z_n)^\gamma\right]\cdot\EE[\mathcal T_{n, \gamma}]<\infty.
\]

Therefore
the series of $\EE\left[(Z_n)^\gamma\right]$ also converges, 
\begin{equation}\label{eq:*}
\sum_{n\geq1}\EE\left[(Z_n)^\gamma\right]<\infty.
\end{equation}
Further, for $\gamma\leq1$ we have 
\[\begin{aligned}
\sum_{n\geq1}\EE\left[\left(Z_n\cdot\mathcal T_{n, 1}\right)\land1\right]\leq&
\sum_{n\geq1}\EE\left[(Z_n\cdot\mathcal T_{n, 1})^\gamma\land1\right]\\
\leq&\sum_{n\geq1}\EE\left[(Z_n\cdot\mathcal T_{n, 1})^\gamma\right]
=\sum_{n\geq1}\EE\left[(Z_n)^\gamma\right]\cdot\EE\left[(\mathcal T_{n, 1})^\gamma\right].
\end{aligned}\]
Moreover,
\[\begin{aligned}
0<\EE\left[\left(\mathcal T_{n,1}\right)^\gamma\right]
=&\EE\left[\left(\int_0^\infty\la_{\ep_n}\rme^{-\la_{\ep_n}t}\rmd t
\int_0^t\rme^{-\eta_{\ep_n}(s)}\rmd s\right)^\gamma\right]\\
=&\EE\left[\left(\int_0^\infty\rme^{-\eta_{\ep_n}(s)-\la_{\ep_n}s}\rmd s\right)^\gamma\right].
\end{aligned}\]
By H\"older's inequality, for $\gamma<1$
\[\begin{aligned}
\EE\left[\left(\mathcal T_{n,1}\right)^\gamma\right]=&
\EE\left[\left(\int_0^\infty\rme^{-\eta_{\ep_n}(s)-\la_{\ep_n}s}\rmd s\right)^\gamma\right]\\
\leq&\left(\EE\int_0^\infty\rme^{-\eta_{\ep_n}(s)-\la_{\ep_n}s}\rmd s\right)^\gamma
=\left(\frac{1}{\la_{\ep_n}+\ell_{\ep_n}(1)}\right)^\gamma.
\end{aligned}\]
By \eqref{eq:*} the   series 
\[
\sum_{n\geq1}\EE\left[\left(Z_n\cdot\mathcal T_{n, 1}\right)\land1\right]\leq\sum_{n\geq1}\EE\left[(Z_n)^\gamma\right]\cdot\EE\left[(\mathcal T_{n, 1})^\gamma\right]\]
converges,
since $\EE\left[(\mathcal T_{n, 1})^\gamma\right]<\infty$ are positive and alternating.

 Therefore,  by \cite[Proposition 3.14]{Kallenberg}
\[
I_\infty=\mathcal T_{0,1}+\sum_{n\geq1}Z_n\cdot\mathcal T_{n,1}<\infty,\qquad a. s. 
\]

\begin{lem} \label{lem3}
$\EE \left[I_\infty(\gamma)\right]<\infty$  if and only if
\begin{align}\label{eq:Tr}
\la_0+\la_1+\ell_0(\gamma)+\ell_1(\gamma)&>0,\\
\ell_0(\gamma)\ell_1(\gamma)+\la_0\ell_1(\gamma)+\la_1\ell_0(\gamma)
+\la_0\la_1\left[1-\widetilde h_0(\gamma)\widetilde h_1(\gamma)\right]&>0.\label{eq:Det}
\end{align}
\end{lem}
\proof
By   \eqref{eq2:Lexp}
$\EE\left[\exp\left(-\gamma X(t)\right)\right]\to0$ as $t\to\infty,$ and 
\[
\EE\left[\int_0^\infty\rme^{-\gamma X(t)}\rmd t\right]<\infty,
\]
if and only if  both eigenvalues $\alpha_1,\;\alpha_2$ of the matrix $\mathcal L(\gamma)$
are positive, i. e. in the case $\mathrm{Tr}\mathcal L(\gamma)>0$ and $\mathrm{Det}\mathcal L(\gamma)>0.$
 This is equivalent to \eqref{eq:Tr}-\eqref{eq:Det}.

Conversely, if  conditions \eqref{eq:Tr}-\eqref{eq:Det} are met, then by \eqref{eq2:Lexp}
the expectation
\[
\EE\left[I_\infty(\gamma)~|~\ep(0)=i\right]
=\int_0^\infty\EE\left[\rme^{-\gamma X(t)}~|~\ep(0)=i\right]\rmd t
=\alpha_1^{-1}{e_1}^{(i)}+\alpha_2^{-1}{e_2}^{(i)},
\]
is   finite, $i\in\{0, 1\}.$ 

Note that if $\gamma $ satisfies condition \eqref{eq:l01gamma01}, then 
condition \eqref{eq:Tr} holds. 
We have the following result.

\begin{theo} \label{theo:++}
Let
\eqref{eq:Det} be  true for some  $\gamma,\;\gamma\in(0, 1],$
satisfying \eqref{eq:l01gamma01}.

Therefore the exponential functional $I_\infty=I_\infty(1)$ is a.s. finite.
\end{theo}

We consider the following tractable examples.
Let the  distributions of jumps $\{Y_n\}$ satisfy
\begin{equation}
\label{eq:Y-beta}
\EE\left[\rme^{-\xi (Y_0+Y_1)}\right]\sim1-b\xi^\beta,\qquad \xi\to0,\qquad \beta>0.
\end{equation}


\begin{ex}
Let $X$ is based on  independent Poisson processes $\eta_n(t)=N_n(t)$
with alternating parameters $\mu_0,\;\mu_1>0,$ that is
their L\'evy-Laplace exponents are given by $\ell_i(\gamma)=\mu_i\left(1-\rme^{-\gamma}\right)\sim\mu_i\gamma,\;\gamma\to0.$
Condition \eqref{eq:l01gamma01} holds for $\gamma>0$.

Let jump magnitudes follow \eqref{eq:Y-beta}
and  $\beta>1.$
In this case for sufficiently small $\gamma>0$ condition \eqref{eq:Det} becomes
\begin{equation}
\label{eq:Det-ex1}
\mu_0\mu_1\gamma^2+(\la_0\mu_1+\la_1\mu_0)\gamma+\la_0\la_1b\gamma^\beta>0
\end{equation}
Note that
 inequality 
\eqref{eq:Det-ex1} has a solution $\gamma,\;\gamma\in(0, 1],$  since
    $\la_0\mu_1+\la_1\mu_0>0.$ Thus $I_\infty<\infty$ a. s.
   
 If  \eqref{eq:Y-beta} holds with 
 $\beta=1,$ then for  \eqref{eq:Det-ex1} it is sufficient to assume that 
 \[
 b\geq-\left(\frac{\mu_0}{\la_0}+\frac{\mu_1}{\la_1}\right).
 \]
  If $\beta<1$ and $b>0,$ then the solution $\gamma,\;0<\gamma\ll1$ of   \eqref{eq:Det-ex1} also exists.
  
  Finally, notice that if $Y_0+Y_1=0$ a. s\textup{., } that is $b=0,$ then \eqref{eq:Det-ex1} holds\textup{,} so
  $I_\infty<\infty$  a. s.
\end{ex}

\begin{ex}\label{ex:telegraph}
Consider the integrated telegraph process $X$  with the states $(\la_i, c_i),$ $i\in\{0, 1\}$,
accompanying with jumps $\{Y_n\}$ which satisfy \eqref{eq:Y-beta}.
The L\'evy-Laplace exponents are 
$\ell_0(\gamma)=c_0\gamma,\;\ell_1(\gamma)=c_1\gamma$ and 
 for sufficiently small $\gamma$
condition \eqref{eq:Det} becomes
\[
c_0c_1\gamma^2+(c_0\la_1+c_1\la_0)\gamma+\la_0\la_1b\gamma^\beta>0.
\]
The exponential functional $I_\infty$ is a. s. finite
in the following three cases\textup{:}
\begin{enumerate}
  \item if $\beta>1$ and $c_0\la_1+c_1\la_0>0;$
  \item if $\beta=1$ and $b\geq-\left(\dfrac{c_0}{\la_0}+\dfrac{c_1}{\la_1}\right);$
    \item if $\beta<1$ and $b>0.$ 
\end{enumerate}

If $Y_0+Y_1=0$ a. s. (that is $b=0$), then     
condition \eqref{eq:Det}  holds  if 
\begin{equation}
\label{eq:telegraph-exponential}
c_0\la_1+c_1\la_0>0.
\end{equation}
\end{ex}

The next example is a bit more sophisticated.

Let $a_0, a_1>0,\; \alpha_0, \alpha_1\in(0, 1),\;\alpha_0\geq\alpha_1$.
Let $\{\eta_n\}$ be independent  stable subordinators with the alternating
 L\'evy-Laplace exponents
$\ell_0(\xi)=a_0\xi^{\alpha_0}$ and $\ell_1(\xi)=a_1\xi^{\alpha_1},\;\xi\geq0.$

Let $X^+$ be the Markov-modulated L\'evy process  
based on such L\'evy processes $\{\eta_n\},$ \eqref{def:X-jump}\textup{,}
and 
$X^-_0$ and $X^-_1$ be the Markov-modulated L\'evy processes  
based on   $\pm\eta_n$ with alternating signs\textup{,} such 
that  the L\'evy-Laplace exponents
of the underlying L\'evy blocks  are
$-\ell_0(\xi)$ and $\ell_1(\xi)$ and $\ell_0(\xi)$ and $-\ell_1(\xi)$ respectively.

Since  \eqref{eq:l01gamma01}  holds for $\pm\ell$ \textup{(}for sufficiently small $\gamma,\;\gamma>0),$
we have the following result.

\begin{theo}
Let $I^+_\infty$ and $I^-_{0,\infty},\;I^-_{1,\infty}$  be  the exponential functionals \eqref{def:I}
of $X^+,\;$  $X^-_0$  and  $X^-_1$ respectively. 
\begin{enumerate}
\item
$I^+_\infty<\infty$ a.s. in 
the following cases\textup{:}
\begin{enumerate}
  \item $\beta>
  \alpha_1;$
  \item $\beta=
  \alpha_1$ and 
  $a_1+\la_1b\geq0,$ if $\alpha_0>\alpha_1,$ or 
  $\la_0a_1+\la_1a_0+\la_0\la_1b\geq0,$ if $\alpha_0=\alpha_1;$
  \item $\beta<
  \alpha_1$ and $b\geq0$.
\end{enumerate}
\item
$I^-_{0, \infty}<\infty$ a.s. in 
the following cases\textup{:}
\begin{enumerate}
  \item $\beta>\alpha_1;$
  \item $\beta=\alpha_1$ and $a_1+\la_1b>0;$
  \item $\beta<\alpha_1$ and $b\geq0$.
  \end{enumerate}
  \item
$I^-_{1, \infty}<\infty$ a.s. in 
the following cases\textup{:}
\begin{enumerate}
  \item $\beta=\alpha_1$ and $\la_1b-a_1\geq0;$
  \item $\beta<\alpha_1$ and $b\geq0$.
\end{enumerate}
\end{enumerate}
\end{theo}

\proof
%
Consider first $I^+_\infty.$ By Theorem \ref{theo:++}
we need to check
\begin{equation}
    a_0a_1\gamma^{\alpha_0+\alpha_1}+\la_0a_1\gamma^{\alpha_1}
    +\la_1a_0\gamma^{\alpha_0}+\la_0\la_1\left[1-\EE\left(\rme^{-\gamma(Y_0+Y_1)}\right)\right]>0, 
    \label{eq:Det-gamma} 
\end{equation}
with some $\gamma,\;\gamma\in(0,\; 1].$
This is equivalent to
\begin{equation}\label{eq:Det-gamma-lim} 
a_0a_1+\frac{\la_0a_1}{\gamma^{\alpha_0}}+\frac{\la_1a_0}{\gamma^{\alpha_1}}
+\frac{\la_0\la_1b}{\gamma^{\alpha_0+\alpha_1-\beta}}+o(1)>0,
\end{equation}
as $\gamma\to0.$ It is easy to see that 
in the cases (1a),\;(1b) and (1c) inequality
\eqref{eq:Det-gamma-lim} holds 
for sufficiently small
$\gamma,\;\gamma>0.$

Second, $I^-_{0, \infty}<\infty$ a. s., if
for sufficiently small $\gamma$
\begin{equation*}
 -a_0a_1-\frac{\la_1a_0}{\gamma^{\alpha_1}}+\frac{\la_0a_1}{\gamma^{\alpha_0}}
+\frac{\la_0\la_1b}{\gamma^{\alpha_0+\alpha_1-\beta}}>0,
\end{equation*}
which is true in the cases (2a),\;(2b) and (2c).

In the case of $I^-_{1,\infty}$ \eqref{eq:Det-gamma}  becomes
\begin{equation}
\label{eq:det-gamma11}
 -a_0a_1-\frac{\la_0a_1}{\gamma^{\alpha_0}}+\frac{\la_1a_0}{\gamma^{\alpha_1}}
+\frac{\la_0\la_1b}{\gamma^{\alpha_0+\alpha_1-\beta}}>0,
\end{equation}
which holds in the cases of (3a) and (3b).
Note that
 if $\beta>\alpha_1$ inequality \eqref{eq:det-gamma11}
does not asymptotically hold (as $\gamma\to0$). 

We will denote by $I^{(0)}_\infty$ and $I^{(1)}_\infty$ the functionals \eqref{def:I},
when the alternating L\'evy process  $X$ starts from the state 0 and 1 respectively.
Let $T$ be the first switching time.
By \eqref{def:X-jump} we have the following identities in law:
\begin{align}
I^{(0)}_\infty\stackrel{D}{=}&\mathcal T^{(0)}(T)+\exp(-Y_0-\eta_0(T))\widetilde I^{(1)}_\infty,\nonumber \\
I^{(1)}_\infty\stackrel{D}{=}&\mathcal T^{(1)}(T)+\exp(-Y_1-\eta_1(T))\widetilde I^{(0)}_\infty,\label{eq:I-jump}
\end{align}
where 
\[
\mathcal T^{(0)}(T)=\int_0^T\rme^{-\eta_0(t)}\rmd t,\qquad
\mathcal T^{(1)}(T)=\int_0^T\rme^{-\eta_1(t)}\rmd t,
\]
 and 
$\widetilde I^{(0)}_\infty$ and $\widetilde I^{(1)}_\infty$ are independent copies of 
$I^{(0)}_\infty$ and $I^{(1)}_\infty$. By identities \eqref{eq:I-jump} for any $t,\;t>0,$  
\begin{equation*}
\left\{
\begin{aligned}
   \PP\{I_\infty^{(0)}>t\} &=\PP\{I_\infty^{(1)}>\exp(Y_0+\eta_0(T))\left(t-\mathcal T^{(0)}(T)\right)\},   \\
   \PP\{I_\infty^{(1)}>t\} &=\PP\{I_\infty^{(0)}>\exp(Y_1+\eta_1(T))\left(t-\mathcal T^{(1)}(T)\right)\}.  
\end{aligned}\right.
\end{equation*}

Since $T$ is exponentially distributed, 
 the density functions $f_0(t)$ and $f_1(t),$ $t\geq0,$ 
 \[
f_0(t)\rmd t=\PP\{I_\infty \in\rmd t~|~\ep(0)=0\},\qquad
f_1(t)\rmd t=\PP\{I_\infty \in\rmd t~|~\ep(0)=1\},
\]
of the distributions 
of $I_\infty^{(0)}$ and $I_\infty^{(1)}$ follow the coupled integral equations:
\begin{equation}
\label{eq:intI-jump}
\left\{
\begin{aligned}
    f_0(t)&=\int_0^\infty\la_0\rme^{-\la_0\tau}
    \EE\left[\rme^{Y_0+\eta_0(\tau)}
    f_1\left(\rme^{Y_0+\eta_0(\tau)}(t-\mathcal T^{(0)}(\tau))\right)\right]\rmd\tau,   \\
    f_1(t)&=\int_0^\infty\la_1\rme^{-\la_1\tau}
    \EE\left[\rme^{Y_1+\eta_1(\tau)}
    f_0\left(\rme^{Y_1+\eta_1(\tau)}(t-\mathcal T^{(1)}(\tau))\right)\right]\rmd\tau. 
\end{aligned}
\right.
\end{equation}

In some particular cases the distributions of $I_\infty^{(0)}$ and $I_\infty^{(1)}$ can be written explicitly by solving system \eqref{eq:intI-jump}. For instance, system \eqref{eq:intI-jump} can be solved in the important particular case 
 of the integrated telegraph process $X$ 
with the  states $(\la_i, c_i),\;i\in\{0, 1\},\;c_0>c_1,$ 
accompanying with the deterministic jumps $y_0, y_1.$

\subsection{Exponential functional for the jump-telegraph process}\label{sec:exponential2}
Consider  an example of the Kac-L\'evy process with jumps and 
with the L\'evy-Laplace exponents $\ell_0(\xi)=c_0\xi,\;\ell_1(\xi)=c_1\xi.$  
Note that if $c_0,\; c_1\leq0,$ then $I^{(0)}_\infty=\infty,\; I^{(1)}_\infty=\infty,\; a.s.$
On the contrary, if both $c_0,\; c_1 $ are positive and $y_0+y_1\geq0,$ 
then the variables $I_\infty^{(0)}$ and $I_\infty^{(1)}$ are a.s. bounded.

In what follows we assume $c_0>c_1,\;c_0>0$ and 
\begin{equation}\label{*}
y_0+y_1=0,\qquad \;y_0\leq0.
\end{equation}
Note that by \eqref{*}
the following inequalities hold   $\forall t,\;t>0,$ a.s.
\begin{align}
\label{ineq:X0}
    y_0+c_1t<&X^{(0)}(t)<c_0t,   \\
    c_1t<&X^{(1)}(t)<y_1+c_0t,
    \label{ineq:X1}  
\end{align}
and 
system \eqref{eq:intI-jump} becomes
\begin{equation}
\left\{
\label{eq:int1-0}
\begin{aligned}
   f_0(t) &=\int_0^\infty\la_0\rme^{-\la_0\tau+y_0+c_0\tau}
   f_1\left(\rme^{y_0+c_0\tau}(t-\mathcal T^{(0)}(\tau))\right)\rmd\tau,\\
   f_1(t) &=\int_0^\infty\la_1\rme^{-\la_1\tau+y_1+c_1\tau}
   f_0\left(\rme^{y_1+c_1\tau}(t-\mathcal T^{(1)}(\tau))\right)\rmd\tau.
\end{aligned}\right.\qquad t>0.
\end{equation}
Here 
\[
\mathcal T^{(0)}(\tau)=\frac{1-\rme^{-c_0\tau}}{c_0},\qquad 
\mathcal T^{(1)}(\tau)=\frac{1-\rme^{-c_1\tau}}{c_1}
\]
(if $c_1=0,$ then $\mathcal T^{(1)}(\tau)=\tau$). 

First, we present the formulae for the density functions $f_0$ and $f_1$
in the case of nonnegative trends, $c_0>c_1\geq0$.
\begin{theo}\label{theo:exp-telegraph-nonnegative}\textbf{Nonnegative $c_0$ and $c_1$.}
Assume that
$
c_0>c_1\geq0.
$ 

Therefore, $I^{(0)}_\infty<\infty,\; I^{(1)}_\infty<\infty,\; a.s.$

\begin{itemize}
\item
If $0<c_1<c_0,$  then
the density functions $f_0$ and $f_1$ have a compact support\textup{,}
and they are given by
\begin{align}
f_0(t)=&A_0(t-a)^{\alpha-1}(b\rme^{-y_0}-t)^{\beta}
\1_{\{a<t<b\rme^{-y_0}\}},\nonumber
\\
f_1(t)=&A_1(t-a\rme^{-y_1})^{\alpha}(b-t)^{\beta-1}
\1_{\{a\rme^{-y_1}<t<b  \}},\label{eq:f0f1+}
\end{align}
where $a=1/c_0,\;\alpha=\la_0/c_0,\;b=1/c_1, \beta=\la_1/c_1$ and
\[A_0=\dfrac{(b\rme^{-y_0}-a)^{-\alpha-\beta}}{\mathrm B(\alpha, \beta+1)},\qquad 
A_1=\dfrac{(b-a\rme^{-y_1})^{-\alpha-\beta}}{\mathrm B(\alpha+1, \beta)};\]
$\mathrm B(\cdot, \cdot)$ is  beta-function.

\item
If $0=c_1<c_0,$ then
the density functions $f_0$ and $f_1$ 
are given by
\begin{align}
    f_0(t)=&\frac{(\la_1\rme^{y_0})^\alpha}{\Gamma(\alpha)}
    \left(t-a\right)^{\alpha-1}\exp(-\la_1\rme^{y_0}(t-a))\1_{\{t>a\}}, \nonumber  \\
    f_1(t)=&\frac{\la_1^{\alpha+1}}{\Gamma(\alpha+1)}
    (t-a\rme^{-y_1})^\alpha\exp\left(-\la_1(t-a\rme^{-y_1})\right)\1_{\{t>a\rme^{-y_1}\}},\label{eq:f0f1+0}
\end{align}
where $a=1/c_0,\;\alpha=\la_0/c_0$ and $\Gamma(\cdot)$ is gamma-function.
\end{itemize}
\end{theo}

\begin{figure}[ht]
\begin{center}
  \includegraphics[scale=0.66]{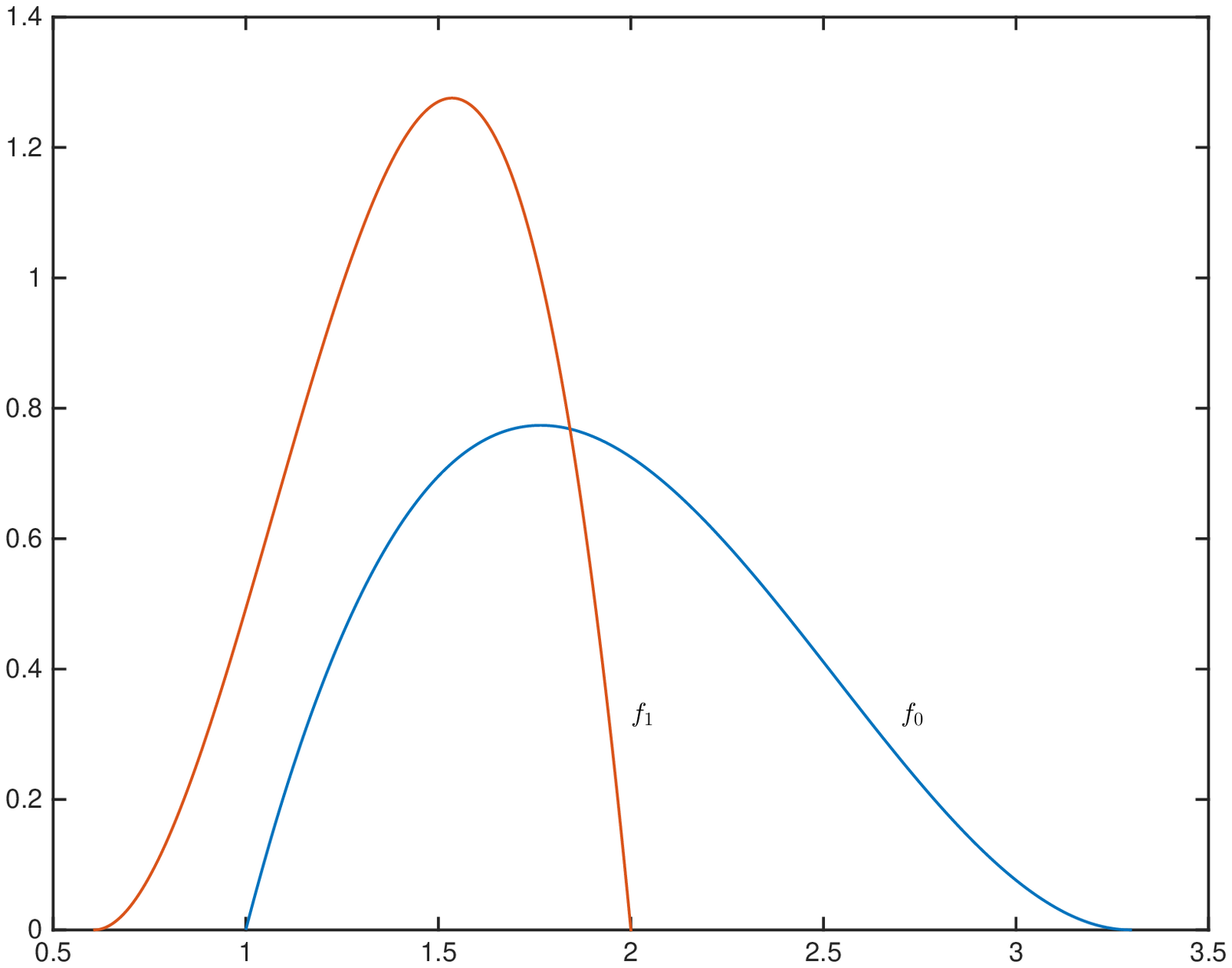}\\
 \caption{Density functions $f_0$ and $f_1,$  \eqref{eq:f0f1+},
 with $\lambda_0=2, \la_1=1;\;c_0=1, c_1=0.5;\; y_0=-0.5, y_1=0.5$.
 {}{}}
 \end{center}
\end{figure}

\proof Let $0<c_1<c_0.$ By \eqref{ineq:X0}-\eqref{ineq:X1}  we have 
\begin{equation*}
    a<I^{(0)}_\infty<b\rme^{-y_0},   \qquad
   a\rme^{-y_1}<I^{(1)}_\infty<b,\qquad a.s.
\end{equation*}
where  $a=1/c_0,\;b=1/c_1.$

By applying the change of variables 
\[\begin{aligned}
u=\rme^{y_0+c_0\tau}\left(t-\mathcal T^{(0)}(\tau)\right)
=&\rme^{y_0}\left[a+\rme^{c_0\tau}(t-a)\right],\quad t\rme^{y_0}<u<b\\
\quad
&\Longleftrightarrow \quad
\tau=a\log\frac{u\rme^{-y_0}-a}{t-a},
\end{aligned}\]
in the first equation of \eqref{eq:int1-0} with $a<t<b\rme^{-y_0}$;
 and 
\[\begin{aligned}
 u=\rme^{y_1+c_1\tau}\left(t-\mathcal T^{(1)}(\tau)\right)
 =&\rme^{y_1}\left[b-\rme^{c_1\tau}(b-t)\right],\quad a<u<t\rme^{y_1}\\
 \quad
&\Longleftrightarrow \quad
\tau=b\log\frac{b-u\rme^{-y_1}}{b-t}
\end{aligned}\]
 in the second one with $a\rme^{-y_1}<t<b,$  
 we found  that system \eqref{eq:int1-0}  is equivalent to 
 \begin{align}
\label{eq:int3-0}
f_0(t)&= \alpha(t-a)^{\alpha-1}
    \int_{t\rme^{y_0}}^{b}(u\rme^{-y_0}-a)^{-\alpha}f_1(u)\rmd u,
    \qquad a<t<b\rme^{-y_0},\\
f_1(t)&=\beta(b-t)^{\beta-1}\int_{a}^{t\rme^{y_1}}(b-u\rme^{-y_1})^{-\beta}f_0(u)\rmd u,
   \qquad a\rme^{-y_1}<t<b.
   \label{eq:int3-1}
\end{align}
 Setting 
 \[
f_0(t)=A_0(t-a)^{\alpha-1}(b-t\rme^{-y_1})^\beta\1_{\{a<t<b\rme^{-y_0}\}}, 
 \]
into \eqref{eq:int3-1} we obtain
\[
f_1(t)=A_1(t-a\rme^{-y_0})^\alpha(b-t)^{\beta-1}\1_{\{a\rme^{-y_1}<t<b\}},
\]
where by \eqref{eq:int3-0} $A_1=\dfrac{\beta A_0}{\alpha}$.
The expressions for $A_0$ and $A_1$  follow  from 
\cite[3.196.3]{GR}.

Let $c_1=0$ and $c_0>0$. Hence, condition \eqref{eq:Tr}-\eqref{eq:Det} 
 (with $\gamma=1$) hold and 
\[
\EE\{I^{(0)}_\infty\}<\infty,\qquad \EE\{I^{(1)}_\infty\}<\infty.
\]
and $I^{(0)}_\infty<\infty,\;I^{(1)}_\infty<\infty,\;a.s.$

System \eqref{eq:int1-0}  
becomes
\begin{equation*}
\left\{
\begin{aligned}
     f_0(t)=&\int_0^\infty\la_0\rme^{-\la_0\tau+y_0+c_0\tau}
    f_1(\rme^{y_0+c_0\tau}(t -\mathcal T^{(0)}(\tau)))\rmd\tau,  
    \qquad a<t,\\
     f_1(t)=&\int_0^\infty\la_1\rme^{-\la_1\tau+y_1}
    f_0(\rme^{y_1}(t-\tau))\rmd\tau,\qquad a\rme^{-y_1}<t,
\end{aligned}
\right.
\end{equation*}
and,  after the corresponding change of variables, 
$\tau=a\log \dfrac{u\rme^{-y-0}-a}{t-a}$ in the first equation and 
$\tau=t-u\rme^{-y_1}$ in the second one,
we get
\begin{equation*}
\left\{
\begin{aligned}
f_0(t)&= \alpha(t-a)^{\alpha-1}
    \int_{t\rme^{y_0}}^{\infty}(u\rme^{-y_0}-a)^{-\alpha}f_1(u)\rmd u,
    \qquad a<t<\infty,\\
f_1(t)&=\la_1\rme^{-\la_1t}\int_a^{t\rme^{y_1}}\exp(\la_1u\rme^{-y_1})f_0(u)\rmd u,
\qquad a\rme^{-y_1}<t.
\end{aligned}
\right.
\end{equation*}

Taking into the second equation
\[
f_0(t)=A_0(t-a)^{\alpha-1}\exp\left(-\la_1t\rme^{y_0}\right)\1_{\{t>a\}}
\]
we easily get
\[
f_1(t)=A_1(t\rme^{y_1}-a)^\alpha\rme^{-\la_1t}\1_{\{t>a\rme^{-y_1}\}},
\]
with $A_1=\la_1A_0/\alpha$.
The expressions 
for $A_0$ and  $A_1$
 follow by the definition of gamma-function:
\[
1=A_0\int_\alpha^\infty(t-a)^{\alpha-1}\rme^{-\la_1t\rme^{y_0}}\rmd t
=A_0\frac{\Gamma(\alpha)\exp\left(-\la_1a\rme^{y_0}\right)}{(\la_1\rme^{y_0})^\alpha}.
\]
The theorem is proved.

The case of $c_0>0>c_1$ with directions of jumps to be   opposite to the sign of current trend
is most important for financial applications, see \cite{KR}. 

\begin{theo}\label{theo:exp-telegraph-negative} 
\textbf{Trends $c_0$ and $c_1$ are of opposite signs.}
 Let $c_1<0<c_0$  and $\alpha=\la_0/c_0>0,$ $\beta=\la_1/c_1<0$. 
 
 \begin{itemize}
\item
If $\alpha+\beta\geq0,$ then 
\begin{equation}
\label{eq:nonconvergence}
\PP\{I^{(0)}_\infty=\infty\}=\PP\{I^{(1)}_\infty=\infty\}=1.
\end{equation}

\item
If $\alpha+\beta<0,$  that is\textup{,}
\begin{equation}
\label{eq:convergence}
\frac{\la_0}{c_0}+\frac{\la_1}{c_1}<0,
\end{equation}
then 
\begin{equation*}
\PP\{I_\infty^{(0)}<\infty\}=\PP\{I_\infty^{(1)}<\infty\}=1.
\end{equation*}

The density functions  
of the distribution in this case are given by
\begin{align}
f_0(t)=&\frac{(a-b\rme^{-y_0})^{-\alpha-\beta}}{B(-\alpha-\beta, \alpha)}
(t-a)^{\alpha-1}(t -b\rme^{-y_0})^\beta\1_{\{t>a\}},\nonumber\\
f_1(t)=&\frac{(a\rme^{-y_1}-b)^{-\alpha-\beta}}{B(-\alpha-\beta, \alpha+1)}
(t-a\rme^{-y_1})^\alpha(t-b)^{\beta-1}\1_{\{t>a\rme^{-y_1}\}}.\label{eq:f0f1pm}
\end{align}
\end{itemize}
\end{theo}

\begin{figure}[ht]
\begin{center}
  \includegraphics[scale=0.66]{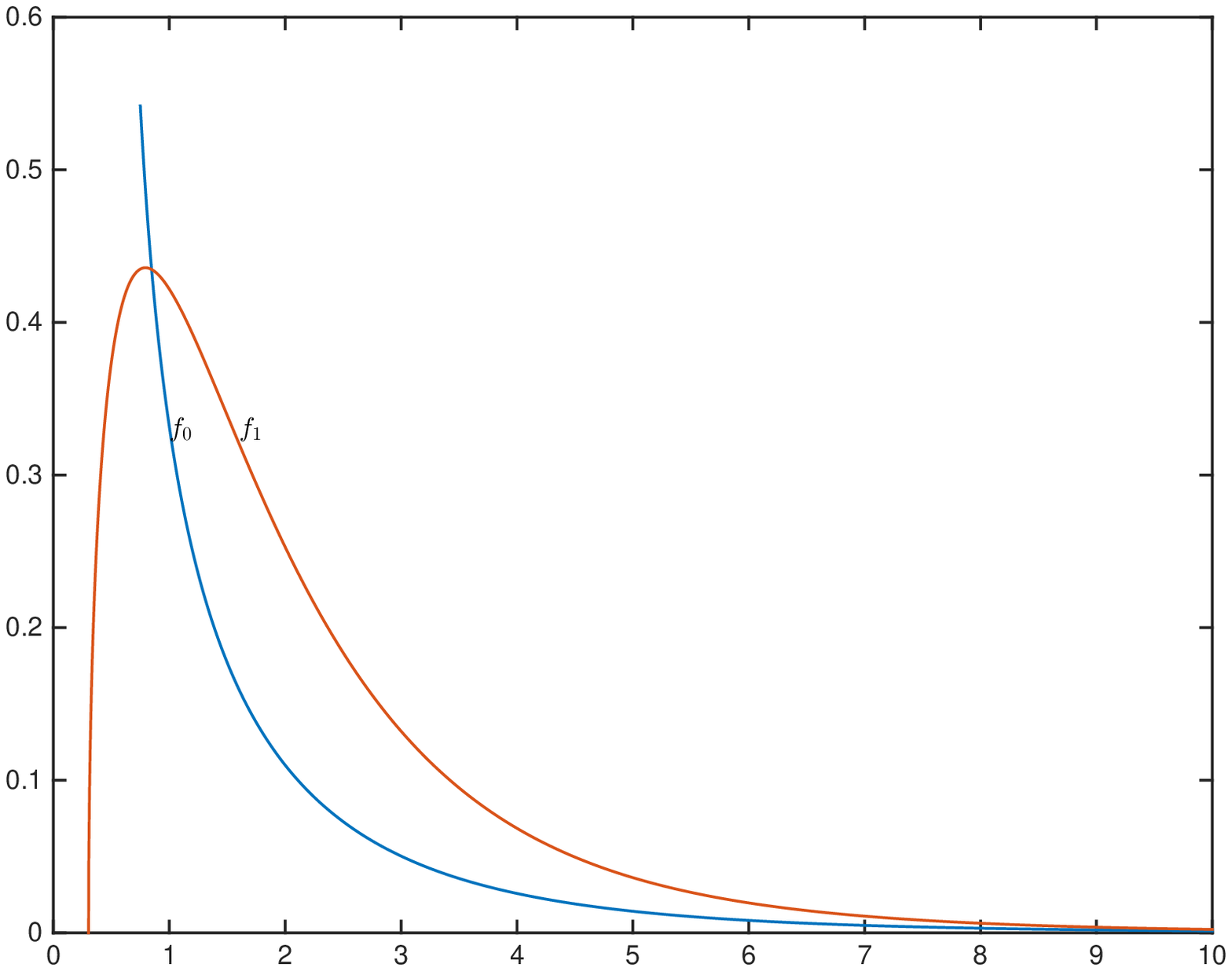}\\
 \caption{Density functions $f_0$ and $f_1$, \eqref{eq:f0f1pm},  
 with $\lambda_0=\la_1=1;\;c_0=2, c_1=-0.1;\; y_0=-0.5, y_1=0.5$.
  {}{}}
 \end{center}
\end{figure}

\proof
System \eqref{eq:int1-0}   for the density functions $f_0$ and $f_1$
after the corresponding  change of variables becomes
\begin{equation}
\label{eq3:intf0f1}
\left\{
\begin{aligned}
f_0(t)&= \alpha(t-a)^{\alpha-1}
    \int_{t\rme^{y_0}}^{\infty}(u\rme^{-y_0}-a)^{-\alpha}f_1(u)\rmd u,
    \qquad a<t,\\
    f_1(t)&=-\beta(t-b)^{\beta-1} 
   \int_{a}^{t\rme^{y_1}} (u\rme^{-y_1}-b)^{-\beta}f_0(u)\rmd u,
   \qquad a\rme^{-y_1}<t.
\end{aligned}
\right.
\end{equation}
Here $\alpha, a>0,\;\beta, b<0.$  

Similarly to the proof of Theorem \ref{theo:exp-telegraph-nonnegative}, 
one can obtain the solution of
\eqref{eq3:intf0f1} in the form
\begin{align}
    f_0(t)=&A(t-a)^{\alpha-1}(t\rme^{-y_1}-b)^\beta\1_{\{t>a\}},  \nonumber \\
    f_1(t)=&\frac{-\beta}{\alpha}A(t\rme^{-y_0}-a)^\alpha(t-b)^{\beta-1} \1_{\{t>a\rme^{-y_1}\}} \label{eq4:intf0f1}
\end{align}
with an indefinite coefficient $A$.

Let  $\alpha+\beta\geq0.$
Note that in this case the integrals 
$\int_\alpha^\infty  f_0(t)\rmd t$ and $\int_{\alpha\rme^{-y_1}}^\infty f_1(t)\rmd t$
of functions $f_0$ and $f_1$ defined by  \eqref{eq:f0f1pm},
diverge, if $A\neq 0$. Hence $A=0$ and we have \eqref{eq:nonconvergence}.

Let $\alpha+\beta<0$. This condition coincides with \eqref{eq:telegraph-exponential}, see
 Example \ref{ex:telegraph}. Thus $I_\infty<\infty$ a. s. and
formulae \eqref{eq:f0f1pm} follow from \eqref{eq4:intf0f1} and  
$
\int_a^\infty f_0(t)\rmd t=1,
$
see \cite[3.196.2]{GR}.

\begin{rem}
The Laplace transforms of these distributions which are obtained in 
Theorem \ref{theo:exp-telegraph-nonnegative} 
and Theorem \ref{theo:exp-telegraph-negative},
\[
\widetilde f_i(s)=\int_0^\infty\rme^{-st}f_i(t)\rmd t,\qquad s>0,
\]
can be expressed in terms
of the  hypergeometric functions\textup{:}
\begin{itemize}
\item if $c_0>c_1>0,$ then
\[
\widetilde f_0(s)=
\negthinspace_1 F_1(\alpha, \alpha+\beta+1;-(b\rme^{-y_0}-a)s)\exp(-as)
\]
and 
\[
\widetilde f_1(s)=
\negthinspace_1 F_1(\alpha+1, \alpha+\beta+1;-(b-a\rme^{-y_1})s)\exp(-a\rme^{-y_1}s);
\]
\item if $c_0>0=c_1,$ then
\[
\widetilde f_0(s)=\left(\frac{\la_1\rme^{y_0}}{\la_1\rme^{y_0}+s}\right)^\alpha\exp(-as)
\]
and 
\[
\widetilde f_1(s)=\left(\frac{\la_1}{\la_1+s}\right)^{\alpha+1}\exp(-a\rme^{-y_1}s);
\]

\item if $c_0>0>c_1,$ then 
\[
\widetilde f_0(s)=\frac{\Gamma(-\beta)}{\Gamma(-\alpha-\beta)}
\Psi(\alpha, \alpha+\beta+1; (a-b\rme^{-y_0})s)\exp(-as)
\]
and 
\[
\widetilde f_1(s)=\frac{\Gamma(-\beta+1)}{\Gamma(-\alpha-\beta)}
\Psi(\alpha+1, \alpha+\beta+1; (a\rme^{-y_1}-b)s)
\exp(-\alpha s\rme^{-y_1}),
\]
where $\beta=\la_1/c_1<0$ and by \eqref{eq:convergence}
$\alpha+\beta=\la_0/c_0+\la_1/c_1<0.$
\end{itemize}
Here $\negthinspace_1 F_1$ is the Kummer
and $\Psi$  the Tricomi  confluent hypergeometric  functions. 
For these formulae  see \textup{\cite[2.1.3]{PBM4},} \eqref{eq:f0f1+}-\eqref{eq:f0f1+0} and
\eqref{eq:f0f1pm}.
\end{rem}




\end{document}